# SEMIPARAMETRICALLY EFFICIENT RANK-BASED INFERENCE FOR SHAPE
## II. OPTIMAL $R$-ESTIMATION OF SHAPE

By Marc Hallin[1], Hannu Oja[2] and Davy Paindaveine[1]

*Université Libre de Bruxelles, University of Tampere and Université Libre de Bruxelles*

A class of $R$-estimators based on the concepts of multivariate signed ranks and the optimal rank-based tests developed in Hallin and Paindaveine [*Ann. Statist.* **34** (2006)] is proposed for the estimation of the shape matrix of an elliptical distribution. These $R$-estimators are root-$n$ consistent under any radial density $g$, without any moment assumptions, and semiparametrically efficient at some prespecified density $f$. When based on normal scores, they are uniformly more efficient than the traditional normal-theory estimator based on empirical covariance matrices (the asymptotic normality of which, moreover, requires finite moments of order four), irrespective of the actual underlying elliptical density. They rely on an original rank-based version of Le Cam's one-step methodology which avoids the unpleasant nonparametric estimation of cross-information quantities that is generally required in the context of $R$-estimation. Although they are not strictly affine-equivariant, they are shown to be equivariant in a weak asymptotic sense. Simulations confirm their feasibility and excellent finite-sample performances.

## 1. Introduction.

1.1. *Rank-based inference for elliptical families.*   An elliptical density over $\mathbb{R}^k$ is determined by a location center $\boldsymbol{\theta} \in \mathbb{R}^k$, a scale parameter $\sigma \in \mathbb{R}_0^+$, a real-valued positive definite symmetric $k \times k$ matrix $\mathbf{V} = (V_{ij})$ with $V_{11} = 1$,

Received November 2004; revised September 2005.
[1]Supported by a P.A.I. contract of the Belgian Federal Government and an Action de Recherche Concertée of the Communauté française de Belgique.
[2]Supported by a grant from the Academy of Finland.

the *shape matrix*, and the so-called *standardized radial density* $g_1$; for a precise definition and comments, see Section 1.2 of [13]. We shall hereafter refer to the latter as HP, further referring to Section HP1.2, Proposition HP2.3, Equation (HP4.5), etc.

Elliptical families have been introduced in multivariate analysis as a reaction against pervasive Gaussian assumptions. Most classical procedures in that field—principal components, discriminant analysis, canonical correlations, multivariate regression, etc.—readily extend to elliptical models, with shape playing the role of covariances or correlations. When $g_1$ is such that the corresponding distribution has finite second-order moments, $\mathbf{V}$ is proportional to the covariance matrix and shape-based procedures coincide with the classical covariance-based ones; unlike covariances, however, shape still makes sense in the absence of moment restrictions. In such a context, robust inference methods, resisting arbitrarily heavy radial tails, are highly desirable and distribution-free rank-based methods naturally come into the picture (see [9, 10, 11, 12] for closely related results).

1.2. *Rank tests.* In the hypothesis-testing context, HP develop a class of semiparametrically optimal signed rank tests for null hypotheses of the form $\mathbf{V} = \mathbf{V}_0$ ($\boldsymbol{\theta}$, $\sigma$ and $g_1$ playing the role of nuisances). Let $\mathbf{X}_1, \ldots, \mathbf{X}_n$ be a random sample from some elliptical distribution characterized by $\boldsymbol{\theta}$, $\sigma$, $\mathbf{V}$ and $g_1$. Assuming that $\boldsymbol{\theta}$ is known (in practice, this $\boldsymbol{\theta}$ can be replaced by any root-$n$ consistent estimate $\hat{\boldsymbol{\theta}}$—see Section HP4.4), denote by $\mathbf{Z}_i := \mathbf{V}_0^{-1/2}(\mathbf{X}_i - \boldsymbol{\theta})$ the $\boldsymbol{\theta}$-centered, $\mathbf{V}_0$-standardized observations. Define the rank $R_i$ as the rank of $d_i := \|\mathbf{Z}_i\|$ among $d_1, \ldots, d_n$ and the multivariate sign $\mathbf{U}_i$ as $\|\mathbf{Z}_i\|^{-1}\mathbf{Z}_i$, $i = 1, \ldots, n$. Considering the matrix-valued signed rank statistic

$$\mathbf{S}_{f_1}(\mathbf{V}_0) := \frac{1}{n}\sum_{i=1}^{n} K_{f_1}\left(\frac{R_i}{n+1}\right)\mathbf{U}_i\mathbf{U}_i',$$

where $K_{f_1} : (0,1) \to \mathbb{R}$ is the score function ensuring optimality at $f_1$, the test statistic developed in HP takes the very simple form [see (HP4.4)]

$$\underset{\sim}{Q}_{f_1}(\mathbf{V}_0) := \frac{nk(k+2)}{2\mathcal{J}_k(f_1)}Q(\mathbf{S}_{f_1}(\mathbf{V}_0)), \qquad \text{where } Q(\mathbf{S}) := \operatorname{tr}(\mathbf{S}^2) - \frac{1}{k}(\operatorname{tr}\mathbf{S})^2.$$
(1.1)

Test procedures based on (1.1) enjoy a number of attractive features: (i) they are valid under arbitrary standardized radial densities $g_1$, irrespective of any moment assumptions, (ii) they are nevertheless (semiparametrically) efficient at some prespecified radial density $f_1$, (iii) they exhibit surprisingly high asymptotic relative efficiencies with respect to classical Gaussian procedures under non-Gaussian $g_1$'s and, quite remarkably, (iv) when Gaussian (van der Waerden) scores are adopted, their ARE's with respect to the classical Gaussian tests [21, 22, 34, 35] are uniformly larger than one; see [38] for this extension of the celebrated Chernoff–Savage [5] result to shape matrices.

These optimality properties, in fact, are all possessed by the noncentrality parameters of the noncentral chi-square asymptotic distributions, under local alternatives, of the rank-based test statistic under consideration. When the radial density, under such alternatives, is $g_1$, these noncentrality parameters are quadratic forms characterized by a symmetric positive definite matrix of the form $\mathcal{J}_k^2(f_1, g_1)\mathcal{J}_k^{-1}(f_1)\mathbf{\Upsilon}_k^{-1}(\mathbf{V})$, where $\mathcal{J}_k(f_1, g_1)$ is a cross-information quantity (cf. (2.7)) and $\mathbf{\Upsilon}_k$ depends on neither $f_1$ nor $g_1$; see Proposition HP4.1. This matrix, for $g_1 = f_1$, coincides with the efficient information matrix $\mathcal{J}_k(f_1)\mathbf{\Upsilon}_k^{-1}(\mathbf{V})$ for $\mathbf{V}$ under $f_1$.

An immediate question which arises is whether such tests have any natural counterparts in the context of point estimation. That is, can we construct estimators $\widehat{\mathbf{V}}^{(n)}$ for the shape matrix that match the performances of those rank-based tests, in the sense of (i) being root-$n$ consistent under any radial density $g_1$, irrespective of any moment assumptions—in sharp contrast with the Gaussian estimators, which require finite second-order moments for consistency and finite fourth-order moments for asymptotic normality, (ii) being nevertheless (semiparametrically) efficient at some prespecified standardized radial density $f_1$ and (iii) exhibiting the same asymptotic relative efficiencies, with respect to classical Gaussian estimators, including (iv) the Chernoff–Savage property of [38]? Such estimators would improve the performance of the existing ones that satisfy the consistency requirement (i), such as Tyler's [45] celebrated affine-equivariant estimator of shape (*scatter*, in Tyler's terminology) $\mathbf{V}_T^{(n)}$ or the estimator of shape based on the *Oja signs* developed in [36]. These estimators are indeed root-$n$ consistent under extremely general conditions (second-order moments, however, are required in [36]), but they are not efficient.

The answer, as we shall see, is positive and the estimators achieving the required performances are $R$-estimators based on the same concepts of multivariate ranks and signs as the test statistics (1.1).

1.3. *R-estimation.* The derivation of such $R$-estimators, however, is by no means straightforward. Traditional $R$-estimators are defined (and computed) via the minimization of some rank-based objective function; see [1, 19, 20, 24, 26] or the review paper by Draper [6]. In the present context, this approach, in connection with (1.1), leads to the definition of an $R$-estimator as

(1.2) $\underset{\sim}{\mathbf{V}}_{f_1}^{(n)} := \underset{\mathbf{V}}{\mathrm{argmin}}\, \underset{\sim}{Q}_{f_1}(\mathbf{V}) = \underset{\mathbf{V}}{\mathrm{argmin}}\left(\mathrm{tr}(\mathbf{S}_{f_1}^2(\mathbf{V})) - \frac{1}{k}(\mathrm{tr}\,\mathbf{S}_{f_1}(\mathbf{V}))^2\right),$

that is, as the value of $\mathbf{V}$ minimizing the sum of squared deviations of the $k$ eigenvalues of the rank-based matrix $\mathbf{S}_{f_1}(\mathbf{V})$ from their arithmetic mean.

This "argmin" definition is intuitively quite appealing. However, from a practical point of view, its implementation is numerically costly when the dimension of the parameter is high [a shape parameter has $k(k+1)/2 - 1$ components]. The same definition is hardly more convenient from a theoretical point of view: as a function of ranks, the objective function $\mathbf{V} \mapsto \underset{\sim}{Q}_{f_1}(\mathbf{V})$ is discontinuous and its monotonicity/convexity properties are all but obvious, so root-$n$ consistency remains a nontrivial issue.

Instead, therefore, we suggest a rank-based adaptation of Le Cam's one-step construction of locally asymptotically optimal estimators. A version, $\underset{\sim}{\boldsymbol{\Delta}}_{f_1}^{(n)}(\mathbf{V})$, measurable with respect to the ranks and signs associated with $\mathbf{V}$, of the semiparametrically efficient (at $\mathbf{V}$ and $f_1$) central sequence for shape can be constructed [see (HP4.1) or (2.6)]; this central sequence is distribution-free with asymptotic covariance matrix $\mathcal{J}_k(f_1) \boldsymbol{\Upsilon}_k^{-1}(\mathbf{V})$. The $f_1$-score version of our $R$-estimator, in vech form (that is, stacking the upper-diagonal elements), is then defined as

$$(1.3) \quad \operatorname{vech}(\underset{\sim}{\mathbf{V}}_{f_1}^{(n)}) := \operatorname{vech}(\mathbf{V}_T^{(n)}) + n^{-1/2}(\alpha^*)^{-1} \left( \boldsymbol{\Upsilon}_k(\mathbf{V}_T^{(n)}) \overset{0}{\underset{\sim}{\boldsymbol{\Delta}}}_{f_1}^{(n)}(\mathbf{V}_T^{(n)}) \right),$$

where $\mathbf{V}_T^{(n)}$ is Tyler's estimator of scatter and $\alpha^*$ is a consistent estimator of the cross-information quantity $\mathcal{J}_k(f_1, g_1)$ [the problem of estimating $\mathcal{J}_k(f_1, g_1)$ is discussed in Section 4]. The resulting $\underset{\sim}{\mathbf{V}}_{f_1}^{(n)}$ is a genuine $R$-estimator since the one-step correction in (1.3) only depends on Tyler's $\mathbf{V}_T^{(n)}$ and the corresponding ranks $R_i$ and signs $\mathbf{U}_i$. Moreover, it is asymptotically equivalent to a random matrix (depending on the actual $g_1$) which is measurable with respect to the ranks and signs associated with the "true" value of $\mathbf{V}$. And if (1.2) admits a root-$n$ consistent sequence of solutions, this sequence of solutions and the one-step definition of $\underset{\sim}{\mathbf{V}}_{f_1}^{(n)}$ are asymptotically equivalent.

The main objective of this paper is to show that $\underset{\sim}{\mathbf{V}}_{f_1}^{(n)}$, as defined in (1.3), indeed satisfies the properties listed under (i)–(iv) which are required of a semiparametrically efficient $R$-estimator.

1.4. *Outline of the paper.* The outline of the paper is as follows. In Section 2, we recall the main definitions related to elliptical symmetry, local asymptotic normality and the relation between ranks and signs on one hand and semiparametric efficiency on the other; whenever possible, we refer to HP for reasons of brevity. Postponing to Section 4 the delicate problem of choosing a consistent estimator $\alpha^*$ for $\mathcal{J}_k(f_1, g_1)$, Section 3 deals with the derivation and asymptotic properties of the one-step $R$-estimator (1.3) based on such arbitrary $\alpha^*$. Section 4 is entirely devoted to the estimation of $\mathcal{J}_k(f_1, g_1)$. We start, in Section 4.1, with a review of the various solutions that have been considered in the literature, explaining why they fail to be fully convincing. Sections 4.2 and 4.3 then propose an original, more sophisticated (yet easily implementable) method inspired by local maximum likelihood ideas. The resulting $R$-estimators enjoy all the asymptotic properties expected from $R$-estimation and, moreover, yield surprisingly high ARE's with respect to the existing methods: see Table 1. These estimators, however, remain unsatisfactory on one count: for fixed sample size $n$, they are not affine-equivariant. They are, nevertheless, equivariant in a *weak asymptotic* sense, as shown in Section 5. A numerical study (Section 6) confirms the excellent performance of the method. The Appendix collects technical proofs.

## 2. Semiparametric efficiency under elliptical symmetry.

2.1. *Uniform local asymptotic normality.* Let $\mathbf{X}^{(n)} := (\mathbf{X}_1^{(n)'}, \ldots, \mathbf{X}_n^{(n)'})'$, $n \in \mathbb{N}$, be a triangular array of $k$-dimensional observations. Let $\mathrm{P}^{(n)}_{\boldsymbol{\theta}, \sigma^2, \mathbf{V}; f_1}$ denote the distribution of $\mathbf{X}^{(n)}$ under the assumption that the $\mathbf{X}_i^{(n)}$'s are i.i.d. with the elliptical density $\underline{f}_{\boldsymbol{\theta}, \sigma^2, \mathbf{V}; f_1}$ described in Section HP1.2 [which we refer to for details, as well as for a precise definition of the parameters $\boldsymbol{\theta}$, $\sigma$, $\mathbf{V}$ and $\boldsymbol{\vartheta}$, the parameter spaces $\boldsymbol{\Theta}$ and $\mathcal{V}_k$, the radial distribution functions $\tilde{F}_1$, the distances $d_i^{(n)}(\boldsymbol{\theta}, \mathbf{V})$, the ranks $R_i^{(n)}(\boldsymbol{\theta}, \mathbf{V})$ and the signs $\mathbf{U}_i^{(n)}(\boldsymbol{\theta}, \mathbf{V})$]. Our objective is the estimation of $\mathbf{V}$ under unspecified $\boldsymbol{\theta}$, $\sigma^2$ and $f_1$.

The relevant statistical experiment involves the nonparametric family

$$(2.1) \quad \mathcal{P}^{(n)} := \bigcup_{f_1 \in \mathcal{F}_A} \mathcal{P}_{f_1}^{(n)} := \bigcup_{f_1 \in \mathcal{F}_A} \{\mathrm{P}^{(n)}_{\boldsymbol{\theta}, \sigma^2, \mathbf{V}; f_1} | \boldsymbol{\theta} \in \mathbb{R}^k, \sigma > 0, \mathbf{V} \in \mathcal{V}_k\},$$

where $f_1$ ranges over the set $\mathcal{F}_A$ of standardized radial densities satisfying Assumptions (A1)–(A2) in HP. The main technical tool is the uniform local asymptotic normality (ULAN), with respect to $\boldsymbol{\vartheta} := (\boldsymbol{\theta}', \sigma^2, (\mathrm{vech}^\circ \mathbf{V})')'$, of the families $\mathcal{P}_{f_1}^{(n)}$. This ULAN property is stated and proved in Section HP2, which we refer to for the definitions of the score functions $\varphi_{f_1}$, $\psi_{f_1}$ and $K_{f_1}$ and for the explicit forms of the central sequences $\boldsymbol{\Delta}_{f_1}^{(n)}(\boldsymbol{\vartheta})$ and information matrices $\boldsymbol{\Gamma}_{f_1}(\boldsymbol{\vartheta})$.

The block-diagonal structure of $\boldsymbol{\Gamma}_{f_1}(\boldsymbol{\vartheta})$ and ULAN imply that substituting (in principle, after adequate discretization) a root-$n$ consistent estimator $\hat{\boldsymbol{\theta}} = \hat{\boldsymbol{\theta}}^{(n)}$ for the unknown location $\boldsymbol{\theta}$ has no influence, asymptotically, on the $\mathbf{V}$-part $\boldsymbol{\Delta}_{f_1;3}^{(n)}$ of the central sequence. Hence, optimal inference about $\mathbf{V}$ can be based, without any loss of (asymptotic) efficiency, on $\boldsymbol{\Delta}_{f_1;3}^{(n)}(\hat{\boldsymbol{\theta}}, \sigma^2, \mathbf{V})$, as if $\hat{\boldsymbol{\theta}}$ were the actual location parameter. This actually follows from the asymptotic linearity property of Section A.1. Therefore, in the derivation of theoretical results, we may tacitly assume, without loss of generality, that $\boldsymbol{\theta} = \mathbf{0}$. The notation $\mathrm{P}_{\sigma^2, \mathbf{V}; f_1}^{(n)}$, $d_i^{(n)}(\mathbf{V})$, $\mathbf{U}_i^{(n)}(\mathbf{V})$, $\boldsymbol{\Delta}_{f_1}^{(n)}(\sigma^2, \mathbf{V})$, $\boldsymbol{\Gamma}_{f_1}(\sigma^2, \mathbf{V})$, etc. will be used in an obvious way instead of $\mathrm{P}_{\mathbf{0}, \sigma^2, \mathbf{V}; f_1}^{(n)}$, $d_i^{(n)}(\mathbf{0}, \mathbf{V})$, $\mathbf{U}_i^{(n)}(\mathbf{0}, \mathbf{V})$, $\boldsymbol{\Delta}_{f_1;3}^{(n)}(\mathbf{0}, \sigma^2, \mathbf{V})$, $\boldsymbol{\Gamma}_{f_1;3}(\mathbf{0}, \sigma^2, \mathbf{V})$, etc. Experiment (2.1) now takes the form

$$(2.2) \quad \mathcal{P}^{(n)} := \bigcup_{f_1 \in \mathcal{F}_A} \mathcal{P}_{f_1}^{(n)} := \bigcup_{f_1 \in \mathcal{F}_A} \bigcup_{\sigma > 0} \mathcal{P}_{\sigma^2; f_1}^{(n)} := \bigcup_{f_1 \in \mathcal{F}_A} \bigcup_{\sigma > 0} \{\mathrm{P}_{\sigma^2, \mathbf{V}; f_1}^{(n)} | \mathbf{V} \in \mathcal{V}_k\}.$$

Although any root-$n$ consistent estimator $\hat{\boldsymbol{\theta}}$ could be used, we suggest adopting the multivariate affine-equivariant median introduced by Hettmansperger and Randles [18] which is itself a "sign-based" estimator. The multivariate signs to be considered, then, are the $\mathbf{U}_i^{(n)}(\hat{\boldsymbol{\theta}}, \mathbf{V})$'s and the ranks to be considered are those of the $d_i^{(n)}(\hat{\boldsymbol{\theta}}, \mathbf{V})$'s.

2.2. *Semiparametric efficiency, ranks and signs.* The partition (2.2) of $\mathcal{P}^{(n)}$ into a collection of parametric subexperiments $\mathcal{P}_{f_1}^{(n)}$, all indexed by $\mathbf{V}$ and $\sigma^2$, induces a semiparametric structure, where $\mathbf{V}$ is the parameter of interest, while $(\sigma^2, f_1)$ plays the role of a nuisance. Except for the unavoidable loss of efficiency resulting from the presence of this nuisance, we would like our estimators to be optimal, that is, to reach semiparametric efficiency bounds, either at some prespecified radial density $f_1$ or at any density belonging to some class $\mathcal{F}_*$ of radial densities.

The semiparametric efficiency bound at $f_1$ is provided by the so-called *efficient information matrix* (see Section HP3.1),

$$\boldsymbol{\Gamma}_{f_1}^*(\mathbf{V}) := \frac{\mathcal{J}_k(f_1)}{4k(k+2)} \mathbf{M}_k (\mathbf{V}^{\otimes 2})^{-1/2} \left[ \mathbf{I}_{k^2} + \mathbf{K}_k - \frac{2}{k} \mathbf{J}_k \right] (\mathbf{V}^{\otimes 2})^{-1/2} \mathbf{M}_k'$$

(2.3)
$$=: \mathcal{J}_k(f_1) \boldsymbol{\Upsilon}_k^{-1}(\mathbf{V});$$

we refer to Section HP1.4 for a definition of the matrices $\mathbf{V}^{\otimes 2}$, $\mathbf{K}_k$, $\mathbf{J}_k$ and $\mathbf{M}_k$, as well as for those of $\mathbf{J}_k^\perp$ and $\mathbf{N}_k$ which we will use later on. This information matrix (2.3) is the asymptotic covariance (under shape matrix $\mathbf{V}$ and density $f_1$) of the *efficient central sequence*

$$(2.4) \quad \boldsymbol{\Delta}_{f_1}^{*(n)}(\mathbf{V}) := \frac{1}{2} n^{-1/2} \mathbf{M}_k (\mathbf{V}^{\otimes 2})^{-1/2} \mathbf{J}_k^\perp \sum_{i=1}^n \varphi_{f_1}\left(\frac{d_i}{\sigma}\right) \frac{d_i}{\sigma} \operatorname{vec}(\mathbf{U}_i \mathbf{U}_i')$$

(see Section HP3.1) which, like $\boldsymbol{\Gamma}^*_{f_1}(\mathbf{V})$, does not depend on $\sigma$ (hence the notation). An estimator $\mathbf{V}^{(n)}$ of $\mathbf{V}$ is semiparametrically efficient at $(\sigma^2, f_1)$ iff the asymptotic distribution under $\mathrm{P}^{(n)}_{\sigma^2,\mathbf{V};f_1}$ of $n^{1/2}\stackrel{\circ}{\mathrm{vech}}(\mathbf{V}^{(n)} - \mathbf{V})$ is the same as that of $(\boldsymbol{\Gamma}^*_{f_1}(\mathbf{V}))^{-1}\boldsymbol{\Delta}^{*(n)}_{f_1}(\mathbf{V})$, that is, iff, under $\mathrm{P}^{(n)}_{\sigma^2,\mathbf{V};f_1}$,

$$(2.5) \qquad n^{1/2}\stackrel{\circ}{\mathrm{vech}}(\mathbf{V}^{(n)} - \mathbf{V}) \stackrel{\mathcal{L}}{\longrightarrow} \mathcal{N}(\mathbf{0}, (\boldsymbol{\Gamma}^*_{f_1}(\mathbf{V}))^{-1}).$$

The difference between $\boldsymbol{\Gamma}_{f_1}(\sigma^2, \mathbf{V})$ and $\boldsymbol{\Gamma}^*_{f_1}(\mathbf{V})$ quantifies the loss of information on $\mathbf{V}$ which is due to the non-specification of $(\sigma^2, f_1)$. It should be emphasized that, whereas this loss depends on the definition of shape (that is, on the arbitrary choice of the normalization $V_{11} = 1$), the semiparametric information bound does not; see Sections HP3.1, HP3.2, [14] and [39] for details.

A general result by Hallin and Werker [17] suggests that, in case

(i) for all $f_1 \in \mathcal{F}_A$ and $\sigma > 0$, the sequence of parametric subexperiments $\mathcal{P}^{(n)}_{\sigma^2;f_1}$ [see (2.2)] is ULAN with central sequence $\boldsymbol{\Delta}^{(n)}_{f_1}(\sigma^2, \mathbf{V})$ and information matrix $\boldsymbol{\Gamma}_{f_1}(\sigma^2, \mathbf{V})$ and

(ii) for all $\mathbf{V} \in \mathcal{V}_k$ and $n \in \mathbb{N}$, the nonparametric subexperiment $\mathcal{P}^{(n)}_{\mathbf{V}} := \{\mathrm{P}^{(n)}_{\sigma^2,\mathbf{V};f_1} | \sigma > 0, f_1 \in \mathcal{F}_A\}$ is generated by a group of transformations $\mathcal{G}^{(n)}_{\mathbf{V}}$ with maximal invariant $\sigma$-field $\mathcal{B}^{(n)}_{\mathbf{V}}$, then the projection $\mathrm{E}[\boldsymbol{\Delta}^{(n)}_{f_1}(\sigma^2, \mathbf{V})| \mathcal{B}^{(n)}_{\mathbf{V}}]$ of $\boldsymbol{\Delta}^{(n)}_{f_1}(\sigma^2, \mathbf{V})$ onto $\mathcal{B}^{(n)}_{\mathbf{V}}$ yields a distribution-free version of the *semiparametrically efficient central sequence* (2.4).

In the present context, this double structure exists: condition (i) is an immediate consequence of Proposition HP2.1 and the generating groups $\mathcal{G}^{(n)}_{\mathbf{V}}$ are the groups of order-preserving radial transformations described in Section HP4.1, which admit the ranks $R_i = R^{(n)}_i(\mathbf{V})$ of the distances $d^{(n)}_i(\mathbf{V})$ and the multivariate signs $\mathbf{U}_i = \mathbf{U}^{(n)}_i(\mathbf{V})$ as maximal invariants. Moreover, $\mathrm{E}[\boldsymbol{\Delta}^{(n)}_{f_1}(\sigma^2, \mathbf{V})| R_1, \ldots, R_n, \mathbf{U}_1, \ldots, \mathbf{U}_n]$ is asymptotically equivalent to

$$\underset{\sim}{\boldsymbol{\Delta}}^{(n)}_{f_1}(\mathbf{V}) := \frac{1}{2}n^{-1/2}\mathbf{M}_k(\mathbf{V}^{\otimes 2})^{-1/2}\mathbf{J}^{\perp}_k \sum_{i=1}^{n} K_{f_1}\left(\frac{R_i}{n+1}\right)\mathrm{vec}(\mathbf{U}_i\mathbf{U}'_i)$$

$$(2.6)$$

$$= \frac{1}{2}n^{-1/2}\mathbf{M}_k(\mathbf{V}^{\otimes 2})^{-1/2}\sum_{i=1}^{n}\left[K_{f_1}\left(\frac{R_i}{n+1}\right)\mathrm{vec}(\mathbf{U}_i\mathbf{U}'_i) - \frac{m^{(n)}_{f_1}}{k}\mathrm{vec}(\mathbf{I}_k)\right]$$

(see Lemma HP4.1), with $K_{f_1}(u) := \varphi_{f_1}(\tilde{F}^{-1}_1(u))\tilde{F}^{-1}_1(u)$ and exact centerings $m^{(n)}_{f_1} := \frac{1}{n}\sum_{i=1}^{n} K_{f_1}(i/(n+1))$.

The properties of $\underset{\sim}{\boldsymbol{\Delta}}_{f_1}^{(n)}(\mathbf{V})$ are summarized in Proposition 2.1 below. For any $g_1 \in \mathcal{F}_A$, define $\boldsymbol{\Gamma}_{f_1,g_1}^*(\mathbf{V}) := \mathcal{J}_k(f_1, g_1) \boldsymbol{\Upsilon}_k^{-1}(\mathbf{V})$, where

$$(2.7) \qquad \mathcal{J}_k(f_1, g_1) := \int_0^1 K_{f_1}(u) K_{g_1}(u) \, du$$

(a cross-information quantity); the notation $\tilde{G}_{1k}$, $\varphi_{g_1}$ is used in an obvious way. Note that $\mathcal{J}_k(f_1, f_1) = \mathcal{J}_k(f_1)$ so that $\boldsymbol{\Gamma}_{f_1,f_1}^*(\mathbf{V})$ reduces to $\boldsymbol{\Gamma}_{f_1}^*(\mathbf{V})$.

PROPOSITION 2.1. *For any $f \in \mathcal{F}_A$, the rank-based random vector* $\underset{\sim}{\boldsymbol{\Delta}}_{f_1}^{(n)}(\mathbf{V})$

  (i) *is distribution-free under* $\{\mathrm{P}_{\sigma^2,\mathbf{V};g_1}^{(n)} | \sigma > 0, g_1 \in \mathcal{F}\}$, *where $\mathcal{F}$ denotes the class of all possible standardized radial densities;*

  (ii) *is asymptotically equivalent, in* $\mathrm{P}_{\sigma^2,\mathbf{V};g_1}^{(n)}$-*probability for any $g_1 \in \mathcal{F}$, to*

$$(2.8) \quad \boldsymbol{\Delta}_{f_1,g_1}^{*(n)}(\mathbf{V}) := \frac{1}{2} n^{-1/2} \mathbf{M}_k(\mathbf{V}^{\otimes 2})^{-1/2} \mathbf{J}_k^\perp \sum_{i=1}^n K_{f_1}\left(\tilde{G}_{1k}\left(\frac{d_i}{\sigma}\right)\right) \mathrm{vec}(\mathbf{U}_i \mathbf{U}_i'),$$

*hence, in* $\mathrm{P}_{\sigma^2,\mathbf{V};f_1}^{(n)}$-*probability, to the semiparametrically efficient (at $f_1$, for any $\sigma$) central sequence for shape* (2.4);

  (iii) *is asymptotically normal under* $\{\mathrm{P}_{\sigma^2,\mathbf{V};g_1}^{(n)} | \sigma > 0, g_1 \in \mathcal{F}\}$ *with mean zero and covariance matrix* $\boldsymbol{\Gamma}_{f_1}^*(\mathbf{V})$;

  (iv) *is asymptotically normal under* $\mathrm{P}_{\sigma^2,\mathbf{V}+n^{-1/2}\mathbf{v};g_1}^{(n)}$, *with mean* $\boldsymbol{\Gamma}_{f_1,g_1}^*(\mathbf{V}) \times \overset{\circ}{\mathrm{vech}}(\mathbf{v})$ *and covariance matrix* $\boldsymbol{\Gamma}_{f_1}^*(\mathbf{V})$ *for any symmetric matrix $\mathbf{v}$ such that $v_{11} = 0$, any $\sigma > 0$ and any $g_1 \in \mathcal{F}_A$;*

  (v) *satisfies, under* $\mathrm{P}_{\sigma^2,\mathbf{V};g_1}^{(n)}$, *as $n \to \infty$, the asymptotic linearity property*

$$(2.9) \quad \underset{\sim}{\boldsymbol{\Delta}}_{f_1}^{(n)}(\mathbf{V} + n^{-1/2} \mathbf{v}^{(n)}) - \underset{\sim}{\boldsymbol{\Delta}}_{f_1}^{(n)}(\mathbf{V}) = -\boldsymbol{\Gamma}_{f_1,g_1}^*(\mathbf{V}) \overset{\circ}{\mathrm{vech}}(\mathbf{v}^{(n)}) + o_\mathrm{P}(1)$$

*for any bounded sequence $\mathbf{v}^{(n)}$ of symmetric matrices such that $v_{11}^{(n)} = 0$, any $\sigma > 0$ and any $g_1 \in \mathcal{F}_A$.*

PROOF. Part (i): distribution-freeness readily follows from the distribution-freeness, under ellipticity, of the ranks $R_i^{(n)}(\mathbf{V})$ and the signs $\mathbf{U}_i^{(n)}(\mathbf{V})$ with respect to which $\underset{\sim}{\boldsymbol{\Delta}}_{f_1}^{(n)}(\mathbf{V})$ is measurable. Part (ii) is covered by Lemma HP4.1. Parts (iii)–(iv) are established in the proof of Proposition HP4.1. Part (v) follows from the more general result given in Proposition A.1 (see Appendix A.1). □

**3. Optimal one-step $R$-estimation of shape.** Tyler's celebrated estimator of shape, $\mathbf{V}_T^{(n)}$, was introduced by Tyler [45] based on the very simple idea that if $\mathbf{X}$ is elliptical with location $\boldsymbol{\theta}$, then its shape $\mathbf{V}$ is entirely characterized by the fact that $\mathbf{U}(\boldsymbol{\theta}, \mathbf{V}) := \mathbf{V}^{-1/2}(\mathbf{X} - \boldsymbol{\theta})/\|\mathbf{V}^{-1/2}(\mathbf{X} - \boldsymbol{\theta})\|$ is centered, with covariance $(1/k)\mathbf{I}_k$. Accordingly, $\mathbf{V}_T^{(n)}$ is defined as the unique shape matrix satisfying $\frac{1}{n}\sum_{i=1}^n \mathbf{U}_i^{(n)}(\boldsymbol{\theta}, \mathbf{V})(\mathbf{U}_i^{(n)}(\boldsymbol{\theta}, \mathbf{V}))' = \frac{1}{k}\mathbf{I}_k$.

Denote by $\mathbf{V}_\#^{(n)}$ a discretized version of $\mathbf{V}_T^{(n)}$. Such discretizations, which turn root-$n$ consistent preliminary estimators into uniformly root-$n$ consistent ones (see, e.g., Lemma 4.4 in [30] for a typical use), are quite standard in Le Cam's one-step construction of estimators (see [31]), and several of them, characterized by a $\#$ subscript, will appear in the sequel. Denoting by $\lceil x \rceil$ the smallest integer larger than or equal to $x$ and by $c_0$ an arbitrary positive constant that does not depend on $n$, the discretized shape $\mathbf{V}_\#^{(n)}$ can be obtained, for instance, by mapping each entry $v_{ij}^{(n)}/(i,j) \neq (1,1)$ of $\mathbf{V}_T^{(n)}$ onto $v_{ij\#}^{(n)} := c_0^{-1}\text{sign}(v_{ij}^{(n)}) n^{-1/2} \lceil n^{1/2} c_0 |v_{ij}^{(n)}| \rceil$. In practice (where $n = n_0$ is fixed), such discretization is not required (as $c_0$ can be arbitrarily large) and actually makes little sense, as one can always decide to start discretization at $n = n_0 + 1$; see Section 4.3 for practical implementation.

Since $\underset{\sim}{\boldsymbol{\Delta}}_{f_1}^{(n)}(\mathbf{V})$ is a version of the efficient central sequence for shape, Le Cam's classical one-step method suggests estimating $\overset{\circ}{\text{vech}}(\mathbf{V})$ by means of

$$(3.1) \quad \overset{\circ}{\text{vech}}(\underset{\sim}{\mathbf{V}}_{f_1\#}^{(n)}) := \overset{\circ}{\text{vech}}(\mathbf{V}_\#^{(n)}) + n^{-1/2}(\boldsymbol{\Gamma}_{f_1,g_1}^*(\mathbf{V}_\#^{(n)}))^{-1} \underset{\sim}{\boldsymbol{\Delta}}_{f_1}^{(n)}(\mathbf{V}_\#^{(n)}).$$

Such an estimator is semiparametrically efficient at $\mathcal{P}_{f_1}^{(n)}$, in the sense of (2.5). Indeed, in view of Proposition 2.1 and the continuity of $\mathbf{V} \mapsto \boldsymbol{\Gamma}_{f_1,g_1}^*(\mathbf{V})$,

$$n^{1/2}\overset{\circ}{\text{vech}}(\underset{\sim}{\mathbf{V}}_{f_1\#}^{(n)} - \mathbf{V}) = n^{1/2}\overset{\circ}{\text{vech}}(\mathbf{V}_\#^{(n)} - \mathbf{V}) + (\boldsymbol{\Gamma}_{f_1,g_1}^*(\mathbf{V}_\#^{(n)}))^{-1}\underset{\sim}{\boldsymbol{\Delta}}_{f_1}^{(n)}(\mathbf{V}_\#^{(n)})$$

$$= n^{1/2}\overset{\circ}{\text{vech}}(\mathbf{V}_\#^{(n)} - \mathbf{V}) + (\boldsymbol{\Gamma}_{f_1,g_1}^*(\mathbf{V}_\#^{(n)}))^{-1}$$

$$\times (\underset{\sim}{\boldsymbol{\Delta}}_{f_1}^{(n)}(\mathbf{V}) - \boldsymbol{\Gamma}_{f_1,g_1}^*(\mathbf{V}) n^{1/2} \overset{\circ}{\text{vech}}(\mathbf{V}_\#^{(n)} - \mathbf{V})) + o_\text{P}(1)$$

$$(3.2) \qquad\qquad = (\boldsymbol{\Gamma}_{f_1,g_1}^*(\mathbf{V}))^{-1} \underset{\sim}{\boldsymbol{\Delta}}_{f_1}^{(n)}(\mathbf{V}) + o_\text{P}(1)$$

$$(3.3) \qquad\qquad = (\boldsymbol{\Gamma}_{f_1,g_1}^*(\mathbf{V}))^{-1} \boldsymbol{\Delta}_{f_1,g_1}^{*(n)}(\mathbf{V}) + o_\text{P}(1)$$

under $\text{P}_{\sigma^2,\mathbf{V};g_1}^{(n)}$, as $n \to \infty$, where application to $\underset{\sim}{\boldsymbol{\Delta}}_{f_1}^{(n)}(\mathbf{V}_\#^{(n)})$ of the asymptotic linearity property (2.9) is made possible, as usual, by the local discreteness of $\mathbf{V}_\#^{(n)}$. The asymptotic representation (3.3) implies, for $g_1 = f_1$, the efficiency of $\underset{\sim}{\mathbf{V}}_{f_1\#}^{(n)}$, whereas (3.2), by providing for $\underset{\sim}{\mathbf{V}}_{f_1\#}^{(n)}$ an asymptotic representation as a signed-rank-measurable quantity, justifies its status as an $R$-estimator.

A major problem, unfortunately, is that (3.1), via $\boldsymbol{\Gamma}_{f_1,g_1}^*(\mathbf{V}_\#^{(n)})$, involves the unknown cross-information quantity $\mathcal{J}_k(f_1,g_1)$ defined in (2.7); $\underset{\sim}{\mathbf{V}}_{f_1\#}^{(n)}$ is, therefore, just a *pseudo-estimator* which cannot be computed from the observations. In order to obtain a genuine estimator, $\underset{\sim}{\widehat{\mathbf{V}}}_{f_1\#}^{(n)}$, say, a consistent estimator $\alpha^*$ must clearly be substituted for $\mathcal{J}_k(f_1,g_1)$. This estimation of $\mathcal{J}_k(f_1,g_1)$ is absolutely crucial in several respects since it not only explicitly enters the definition of the one-step estimator, but also characterizes its asymptotic covariance. However, obtaining a consistent estimator $\alpha^*$ of $\mathcal{J}_k(f_1,g_1)$—the expectation of a function that depends on the unknown underlying $g_1$—is a delicate problem. Accordingly, we defer the discussion of this issue to Section 4, where, after a review of the various methods available in the literature, we present an original method inspired by local maximum likelihood ideas.

Therefore, in the present section, we define the $f_1$-score $R$-estimator $\underset{\sim}{\widehat{\mathbf{V}}}_{f_1\#}^{(n)}$ as the value of $\underset{\sim}{\mathbf{V}}_{f_1\#}^{(n)}$ resulting from substituting into (3.1) an arbitrary consistent estimator $\alpha^*$ for the unknown $\mathcal{J}_k(f_1,g_1)$. Up to discretization, $\underset{\sim}{\widehat{\mathbf{V}}}_{f_1\#}^{(n)}$ thus is defined as in (1.3). Irrespective of the choice of $\alpha^*$, the resulting one-step $R$-estimators $\underset{\sim}{\widehat{\mathbf{V}}}_{f_1\#}^{(n)}$ are asymptotically equivalent (under $\mathcal{P}^{(n)}$) to the pseudo-estimator $\underset{\sim}{\mathbf{V}}_{f_1\#}^{(n)}$ and, hence, also to the signed rank statistics (3.2) based on the "genuine ranks." Proposition 3.1 summarizes the main properties of these estimators: (i) they are asymptotically equivalent to a function of the genuine ranks and signs, they are asymptotically normal, and their covariance matrix is the inverse of the covariance matrix characterizing the local powers of the optimal rank tests derived in HP, (ii) when based on $f_1$-scores, they are semiparametrically efficient at radial density $f_1$, (iii) for finite $n$, they can be expressed as a linear combination of the Tyler shape matrix and a rank-based shape matrix involving the Tyler ranks and signs, (iv) their asymptotic covariance matrix, under any elliptical density, is proportional to the asymptotic covariance matrices of the Tyler and Gaussian ML estimators. The proportionality constant, which can be considered as a measure of asymptotic relative efficiency, is provided in (v). In order to obtain a simpler expression for the asymptotic covariance matrix of $\mathrm{vec}(\underset{\sim}{\mathbf{V}}_{f_1\#}^{(n)})$ (cf. 3.8), we define $\mathbf{Q}_k(\mathbf{V}) := [k(k+2)]^{-1}\mathbf{M}_k'\boldsymbol{\Upsilon}_k(\mathbf{V})\mathbf{M}_k$. As shown in the proof of Lemma HP3.1 (with $\mathbf{N}_k$ defined in Section HP1.4),

$$\boldsymbol{\Upsilon}_k(\mathbf{V}) = k(k+2)\,\mathbf{N}_k \mathbf{Q}_k(\mathbf{V}) \mathbf{N}_k'. \tag{3.4}$$

PROPOSITION 3.1. *Let $f_1$ and $g_1$ belong to $\mathcal{F}_A$. Then*

(i) *under $\mathrm{P}^{(n)}_{\sigma^2,\mathbf{V};g_1}$, as $n \to \infty$,*

$$n^{1/2}\,\mathring{\mathrm{vech}}(\widehat{\underset{\sim}{\mathbf{V}}}^{(n)}_{f_1\#} - \mathbf{V}) = (\boldsymbol{\Gamma}^*_{f_1,g_1}(\mathbf{V}))^{-1}\underset{\sim}{\boldsymbol{\Delta}}^{(n)}_{f_1}(\mathbf{V}) + o_\mathrm{P}(1) \tag{3.5}$$

$$= (\boldsymbol{\Gamma}^*_{f_1,g_1}(\mathbf{V}))^{-1}\boldsymbol{\Delta}^{*(n)}_{f_1,g_1}(\mathbf{V}) + o_\mathrm{P}(1) \tag{3.6}$$

$$\xrightarrow{\mathcal{L}} \mathcal{N}(\mathbf{0},(\mathcal{J}_k(f_1)/\mathcal{J}_k^2(f_1,g_1))\boldsymbol{\Upsilon}_k(\mathbf{V})) \tag{3.7}$$

*or, in terms of $\mathrm{vec}\,\mathbf{V}$,*

$$n^{1/2}\mathrm{vec}(\widehat{\underset{\sim}{\mathbf{V}}}^{(n)}_{f_1\#} - \mathbf{V}) \xrightarrow{\mathcal{L}} \mathcal{N}(\mathbf{0},(k(k+2)\mathcal{J}_k(f_1)/\mathcal{J}_k^2(f_1,g_1))\mathbf{Q}_k(\mathbf{V})); \tag{3.8}$$

(ii) $\widehat{\underset{\sim}{\mathbf{V}}}^{(n)}_{f_1\#}$ *is semiparametrically efficient at $\{\mathrm{P}^{(n)}_{\sigma^2,\mathbf{V};f_1}\,|\,\sigma>0, \mathbf{V}\in\mathcal{V}_k\}$;*

(iii)

$$\begin{aligned}\widehat{\underset{\sim}{\mathbf{V}}}^{(n)}_{f_1\#} &= \left(1 - \frac{k(k+2)}{\alpha^*}(\underset{\sim}{\mathbf{W}}^{(n)}_{f_1\#})_{11}\right)\mathbf{V}^{(n)}_\# \\ &+ \left(\frac{k(k+2)}{\alpha^*}(\underset{\sim}{\mathbf{W}}^{(n)}_{f_1\#})_{11}\right)\underset{\sim}{\mathbf{W}}^{(n)}_{f_1\#}\big/(\underset{\sim}{\mathbf{W}}^{(n)}_{f_1\#})_{11}\end{aligned} \tag{3.9}$$

*for all $n$, where $\underset{\sim}{\mathbf{W}}^{(n)}_{f_1\#} := \underset{\sim}{\mathbf{W}}^{(n)}_{f_1}(\mathbf{V}^{(n)}_\#)$, with*

$$\underset{\sim}{\mathbf{W}}^{(n)}_{f_1}(\mathbf{V}) := \mathbf{V}^{1/2}\left[\frac{1}{n}\sum_{i=1}^n K_{f_1}\left(\frac{R_i^{(n)}(\mathbf{V})}{n+1}\right)\mathbf{U}^{(n)}_i(\mathbf{V})\mathbf{U}^{(n)'}_i(\mathbf{V})\right]\mathbf{V}^{1/2} \tag{3.10}$$

*and $\alpha^*$ is the consistent estimator of $\mathcal{J}_k(f_1,g_1)$ entering the construction of $\widehat{\underset{\sim}{\mathbf{V}}}^{(n)}_{f_1\#}$;*

(iv) *the Gaussian ML estimator is $\mathbf{V}^{(n)}_\mathcal{G} := \boldsymbol{\Sigma}^{(n)}/(\boldsymbol{\Sigma}^{(n)})_{11}$ with*

$$\boldsymbol{\Sigma}^{(n)} := (n-1)^{-1}\sum_{i=1}^n(\mathbf{X}_i - \bar{\mathbf{X}})(\mathbf{X}_i - \bar{\mathbf{X}})';$$

*provided that the kurtosis coefficient $\kappa_k(g_1) := (kE_k(g_1))/((k+2)D_k^2(g_1)) - 1$ [where we let $E_k(g_1) := \int_0^1(\tilde{G}_{1k}^{-1}(u))^4\,du$ and $D_k(g_1) := \int_0^1(\tilde{G}_{1k}^{-1}(u))^2\,du$] is finite, then under $\mathrm{P}^{(n)}_{\sigma^2,\mathbf{V};g_1}$,*

$$n^{1/2}\mathrm{vec}(\mathbf{V}^{(n)}_\mathcal{G} - \mathbf{V}) \xrightarrow{\mathcal{L}} \mathcal{N}(\mathbf{0},(1+\kappa_k(g_1))\mathbf{Q}_k(\mathbf{V})) \qquad \text{as } n\to\infty;$$

(v) *the ARE (i.e., the inverse ratio of asymptotic variances) under* $\mathrm{P}^{(n)}_{\sigma^2,\mathbf{V};g_1}$, *where* $g_1$ *is such that* $\kappa_k(g_1) < \infty$ *(resp., without any moment assumption on* $g_1$), *of* $\widehat{\underset{\sim}{\mathbf{V}}}^{(n)}_{f_1\#}$ *with respect to* $\mathbf{V}^{(n)}_{\mathcal{G}}$ *(resp., with respect to* $\mathbf{V}^{(n)}_{T}$) *is* $\frac{1+\kappa_k(g_1)}{k(k+2)} \frac{\mathcal{J}^2_k(f_1,g_1)}{\mathcal{J}_k(f_1)}$ *(resp.,* $\frac{1}{k^2} \frac{\mathcal{J}^2_k(f_1,g_1)}{\mathcal{J}_k(f_1)}$).

PROOF. See Appendix (Section A.2). □

Note that the ARE's in part (v) of the proposition are unambiguously defined, despite the multivariate setting, as the asymptotic covariance matrices of (the vec versions of) $\widehat{\underset{\sim}{\mathbf{V}}}^{(n)}_{f_1\#}$, $\mathbf{V}^{(n)}_{\mathcal{G}}$ and $\mathbf{V}^{(n)}_{T}$ all are proportional to $\mathbf{Q}_k(\mathbf{V})$. Their relative performances can thus be described by a single number, a fact that was already observed in [44] (see also [33]); the situation is entirely different for covariance matrices, where two numbers are required [36, 37, 43].

These ARE's coincide with those obtained in HP for the problem of testing $\mathbf{V} = \mathbf{V}_0$ (see Proposition HP4.2). An immediate corollary is that the Chernoff–Savage result of [38] also applies here: the ARE's of the van der Waerden (Gaussian-score) versions $\widehat{\underset{\sim}{\mathbf{V}}}^{(n)}_{\mathrm{vdW}\#}$ of our $R$-estimators ($K_{f_1} = \Psi_k^{-1}$, where $\Psi_k$ stands for the chi-square distribution function with $k$ degrees of freedom—see Section HP4.2) with respect to the Gaussian estimator $\mathbf{V}^{(n)}_{\mathcal{G}}$ are uniformly larger than one (and equal to one only at the multinormal); the Pitman-inadmissibility of $\mathbf{V}^{(n)}_{\mathcal{G}}$ follows.

Table 1 provides some numerical values, under various Student ($t_\nu$) and normal ($\mathcal{N}$) radial densities $g_1$, of the ARE's in Proposition 3.1(v); for details on elliptical Student densities, see Section HP1.2. Note that under Student densities with four degrees of freedom or less, the ARE of $\widehat{\underset{\sim}{\mathbf{V}}}^{(n)}_{f_1\#}$ with respect to $\mathbf{V}^{(n)}_{\mathcal{G}}$ is infinite since $n^{1/2}(\mathbf{V}^{(n)}_{\mathcal{G}} - \mathbf{V})$ is not even $O_{\mathrm{P}}(1)$. Also, note that the limits as $\nu \to 0$ of the ARE's under $t_\nu$, with respect to Tyler's $\mathbf{V}^{(n)}_{T}$, of any $\widehat{\underset{\sim}{\mathbf{V}}}^{(n)}_{\nu_0\#}$ (the $R$-estimator associated with $t_{\nu_0}$ scores) and $\widehat{\underset{\sim}{\mathbf{V}}}^{(n)}_{\mathrm{vdW}\#}$ are relatively modest and strictly less than one; see column $t_0$ in Table 1 for numerical values. In fact,

$$\lim_{\nu \to 0} \mathrm{ARE}_{t_\nu}[\widehat{\underset{\sim}{\mathbf{V}}}^{(n)}_{\nu_0\#} / \mathbf{V}^{(n)}_{T}] = \frac{k(k+\nu_0+2)}{(k+2)(k+\nu_0)} < 1$$

and

$$\lim_{\nu \to 0} \mathrm{ARE}_{t_\nu}[\widehat{\underset{\sim}{\mathbf{V}}}^{(n)}_{\mathrm{vdW}\#} / \mathbf{V}^{(n)}_{T}] = \frac{k}{k+2} < 1.$$

TABLE 1
*ARE's of the rank-based estimators $\widehat{\underset{\sim}{\mathbf{V}}}_{0.5\#}^{(n)}$, $\widehat{\underset{\sim}{\mathbf{V}}}_{3\#}^{(n)}$, $\widehat{\underset{\sim}{\mathbf{V}}}_{10\#}^{(n)}$ and $\widehat{\underset{\sim}{\mathbf{V}}}_{\mathrm{vdW}\#}^{(n)}$ (associated with $t_{0.5}$, $t_3$, $t_{10}$ and Gaussian scores, respectively) with respect to Tyler's $\mathbf{V}_T^{(n)}$ and, in parentheses, with respect to the Gaussian estimator $\mathbf{V}_\mathcal{G}^{(n)}$, under $k$-variate Student densities (with $\nu$ degrees of freedom, $\nu = 0.5, 3, 10$), along with the limiting values obtained for $\nu \to 0$ and $\nu \to \infty$ (the multinormal case), for $k = 2$, $3$, $4$, $6$ and $10$.*

|  | | **Underlying density** | | | | |
|---|---|---|---|---|---|---|
|  | $k$ | $t_0$ | $t_{0.5}$ | $t_3$ | $t_{10}$ | $\mathcal{N}$ |
| $\widehat{\underset{\sim}{\mathbf{V}}}_{0.5\#}^{(n)}$ | 2 | 0.900 ($\infty$) | 1.111 ($\infty$) | 1.246 ($\infty$) | 1.280 (0.853) | 1.296 (0.648) |
|  | 3 | 0.943 ($\infty$) | 1.061 ($\infty$) | 1.145 ($\infty$) | 1.173 (0.939) | 1.189 (0.713) |
|  | 4 | 0.963($\infty$) | 1.038 ($\infty$) | 1.098 ($\infty$) | 1.121 (0.996) | 1.136 (0.757) |
|  | 6 | 0.981 ($\infty$) | 1.020 ($\infty$) | 1.054 ($\infty$) | 1.070 (1.070) | 1.083 (0.813) |
|  | 10 | 0.992 ($\infty$) | 1.008 ($\infty$) | 1.024 ($\infty$) | 1.034 (1.149) | 1.044 (0.870) |
| $\widehat{\underset{\sim}{\mathbf{V}}}_{3\#}^{(n)}$ | 2 | 0.700 ($\infty$) | 0.969 ($\infty$) | 1.429 ($\infty$) | 1.651 (1.101) | 1.792 (0.896) |
|  | 3 | 0.800 ($\infty$) | 0.972 ($\infty$) | 1.250 ($\infty$) | 1.400 (1.120) | 1.507 (0.904) |
|  | 4 | 0.857($\infty$) | 0.977 ($\infty$) | 1.667 ($\infty$) | 1.278 (1.136) | 1.366 (0.911) |
|  | 6 | 0.917 ($\infty$) | 0.985 ($\infty$) | 1.091 ($\infty$) | 1.162 (1.162) | 1.229 (0.921) |
|  | 10 | 0.962 ($\infty$) | 0.992 ($\infty$) | 1.040 ($\infty$) | 1.078 (1.198) | 1.123 (0.936) |
| $\widehat{\underset{\sim}{\mathbf{V}}}_{10\#}^{(n)}$ | 2 | 0.583 ($\infty$) | 0.829 ($\infty$) | 1.376 ($\infty$) | 1.714 (1.143) | 1.961 (0.980) |
|  | 3 | 0.692 ($\infty$) | 0.861 ($\infty$) | 1.212 ($\infty$) | 1.444 (1.156) | 1.633 (0.979) |
|  | 4 | 0.762($\infty$) | 0.887 ($\infty$) | 1.136 ($\infty$) | 1.313 (1.167) | 1.468 (0.979) |
|  | 6 | 0.844 ($\infty$) | 0.921 ($\infty$) | 1.070 ($\infty$) | 1.185 (1.185) | 1.304 (0.978) |
|  | 10 | 0.917 ($\infty$) | 0.955 ($\infty$) | 1.027 ($\infty$) | 1.091 (1.212) | 1.174 (0.978) |
| $\widehat{\underset{\sim}{\mathbf{V}}}_{\mathrm{vdW}\#}^{(n)}$ | 2 | 0.500 ($\infty$) | 0.720 ($\infty$) | 1.280 ($\infty$) | 1.681 (1.120) | 2.000 (1.000) |
|  | 3 | 0.600 ($\infty$) | 0.757 ($\infty$) | 1.130 ($\infty$) | 1.415 (1.132) | 1.667 (1.000) |
|  | 4 | 0.667 ($\infty$) | 0.786 ($\infty$) | 1.063 ($\infty$) | 1.285 (1.142) | 1.500 (1.000) |
|  | 6 | 0.750 ($\infty$) | 0.829 ($\infty$) | 1.005 ($\infty$) | 1.159 (1.159) | 1.333 (1.000) |
|  | 10 | 0.833 ($\infty$) | 0.877 ($\infty$) | 0.973 ($\infty$) | 1.067 (1.186) | 1.200 (1.000) |

This can be explained by the fact that, roughly speaking, "$\mathbf{V}_T^{(n)}$ is optimal at $t_0$." In more rigorous terms, we have that, for any fixed $n$,

$$(3.11) \qquad \widetilde{\underset{\sim}{\mathbf{V}}}_{\nu\#}^{(n)} - \mathbf{V}_T^{(n)} = o(1) \qquad \mathcal{P}^{(n)}\text{-a.s., as } \nu \to 0.$$

Indeed, the scores $K_\nu$ associated with the $k$-dimensional Student $t_\nu$ are $K_\nu(u) = k(k+\nu)G_{k,\nu}^{-1}(u)/(\nu + kG_{k,\nu}^{-1}(u))$, $u \in (0,1)$, where $G_{k,\nu}$ stands for the Fisher–Snedecor distribution function with $k$ and $\nu$ degrees of freedom. It is easily checked that $G_{k,\nu}^{-1}(u)/\nu \to \infty$ as $\nu \to 0$ so that $\lim_{\nu \to 0} K_\nu(u) = k$ for all $u \in (0,1)$. It follows (with obvious notation) that $\underset{\sim}{\mathbf{W}}_{\nu\#}^{(n)} - \mathbf{V}_\#^{(n)} = o(1)$, $\mathcal{P}^{(n)}$ -a.s., as $\nu \to 0$. This, in view of (3.9), implies (3.11). Similarly, it can be shown that (using obvious notation) for all fixed $n$ and $\nu$, $\widehat{\underset{\sim}{\mathbf{V}}}_{\nu\#}^{(n)}(\mathbf{x}_k) -$ $\mathbf{V}_T^{(n)}(\mathbf{x}_k)$ is $o(1)$ as $k \to \infty$ along any sequence $(\mathbf{x}_k, k = 2, 3, \ldots)$, where $\mathbf{x}_k = (\mathbf{x}_{k1}, \ldots, \mathbf{x}_{kn})$ is an $n$-tuple of vectors in $\mathbb{R}^k$; here, for $k > n$, $\mathbf{V}_T^{(n)}(\mathbf{x}_k)$ can be taken as any solution of Tyler's M-equation. This explains the fact that for all fixed $\nu$, the ARE of $\widehat{\underset{\sim}{\mathbf{V}}}_{\nu\#}^{(n)}$ with respect to $\mathbf{V}_T^{(n)}$ goes to 1 as $k \to \infty$. Incidentally, this also holds for the van der Waerden version of our estimators: as the dimension $k$ of the observation space goes to infinity, the information contained in the radii $d_i$ becomes negligible when compared with that contained in the directions $\mathbf{U}_i$.

**4. Estimation of cross-information coefficients.** Our estimators $\widehat{\underset{\sim}{\mathbf{V}}}_{f_1\#}^{(n)}$, thus far, have only been defined up to the choice of a consistent estimator $\alpha^*$ of the unknown cross-information quantity $\mathcal{J}_k(f_1, g_1)$ defined in (2.7). In this section, we first review the various methods available in the literature for estimating $\mathcal{J}_k(f_1, g_1)$ and then present an original method which relies on a local maximum likelihood argument.

4.1. *A brief review of the literature.* The problem of estimating the cross-information coefficient $\mathcal{J}_k(f_1, g_1)$ has always been around in $R$-estimation and probably explains why it has never been as popular as rank tests in applications. Simple consistent estimators of cross-information coefficients (the definition of which depends on the problem under study) have been proposed by Lehmann [32] and Sen [42] for one- and two-sample location problems; these estimators are based on comparisons of confidence interval lengths, a method involving the arbitrary choice of a confidence level $(1-\alpha)$ which has quite an impact on the final result.

Another simple method can be obtained from the asymptotic linearity property of rank statistics (see [2, 29] or [25], page 321 for univariate location and regression). This method extends quite easily to the present context via the asymptotic linearity property (2.9). The latter indeed implies that for all $f_1, g_1 \in \mathcal{F}_A$ and any $k \times k$ symmetric matrix $\mathbf{v}$ such that $v_{11} = 0$,

$$\underset{\sim}{\mathbf{\Delta}}_{f_1}^{(n)}(\mathbf{V}_\#^{(n)} + n^{-1/2}\mathbf{v}) - \underset{\sim}{\mathbf{\Delta}}_{f_1}^{(n)}(\mathbf{V}_\#^{(n)}) = \underset{\sim}{\mathbf{\Delta}}_{f_1}^{(n)}(\mathbf{V} + n^{-1/2}\mathbf{v}) - \underset{\sim}{\mathbf{\Delta}}_{f_1}^{(n)}(\mathbf{V}) + o_P(1)$$

$$= -\mathcal{J}_k(f_1, g_1)\mathbf{\Upsilon}_k^{-1}(\mathbf{V})\overset{\circ}{\text{vech}}(\mathbf{v}) + o_P(1),$$

under $P_{\sigma^2, \mathbf{V}; g_1}^{(n)}$, as $n \to \infty$. Thus, for any $\mathbf{v}$,

$$(4.1) \quad \alpha^*(\mathbf{v}) := \|\underset{\sim}{\mathbf{\Delta}}_{f_1}^{(n)}(\mathbf{V}_\#^{(n)} + n^{-1/2}\mathbf{v}) - \underset{\sim}{\mathbf{\Delta}}_{f_1}^{(n)}(\mathbf{V}_\#^{(n)})\| / \|\mathbf{\Upsilon}_k^{-1}(\mathbf{V}_\#^{(n)})\overset{\circ}{\text{vech}}(\mathbf{v})\|$$

is a consistent estimate, under $P_{\sigma^2, \mathbf{V}; g_1}^{(n)}$, of $\mathcal{J}_k(f_1, g_1)$. This method, however, is likely to suffer the same weaknesses as the univariate traditional idea; in particular, these "naive" estimators involve the arbitrary choice of a "small" perturbation of the parameter [the choice of a particular $\mathbf{v}$ in (4.1) is indeed as good/bad as that of $2\mathbf{v}$, $3\mathbf{v}$, etc.). Theory again provides no guidelines for this choice which, unfortunately, has a dramatic impact on the output.

More elaborate approaches involve a kernel estimate of $g_1$ and, hence, cannot be expected to perform well under small and moderate sample sizes. Such kernel methods have been considered, for Wilcoxon scores, by [41] (see also [3, 4, 7] and, in a more general setting, in Section 4.5 of [27]. They also require arbitrary choices (window width and kernel or, as in [27], the choice of the order $\alpha$ of an empirical quantile) for which universal recommendation seems hardly possible (see [28] for an empirical investigation). Moreover, estimating the actual underlying density is somewhat incompatible with the group-invariance spirit of the rank-based approach: if, indeed, the unknown density $g_1$ is eventually to be estimated by some $\hat{g}_1$, then why not simply adopt a more traditional estimated-score approach based on the asymptotic reconstruction, via $\boldsymbol{\Delta}_{\hat{g}_1}^{*(n)}$, of the efficient central sequence $\boldsymbol{\Delta}_{g_1}^{*(n)}$?

4.2. *An original* (*local likelihood*) *method*: *consistency and efficiency.* A more sophisticated way of dealing with the estimation of $\mathcal{J}_k(f_1, g_1)$ can be obtained by further exploiting the ULAN structure of the model. The basic intuition is that of solving a local likelihood equation. Consistency, however, requires somewhat confusing discretization steps which, as usual, are needed in formal proofs only. Therefore, we provide two descriptions of the method: this section carefully covers the details of discretization and establishes the consistency of the proposed estimator (hence, that of the resulting $\widehat{\underset{\sim}{\mathbf{V}}}_{f_1\#}^{(n)}$), while Section 4.3 below, where discretization is skipped, can be used for practical implementation.

Consider the sequence of (random) half-lines,

$$\mathcal{D}_{\#}^{(n)} = \mathcal{D}_{\#}^{(n)}(\mathbf{V}_{\#}^{(n)}; \underset{\sim}{\boldsymbol{\Delta}}_{f_1}^{(n)}(\mathbf{V}_{\#}^{(n)})) = \{\mathring{\text{vech}}(\underset{\sim}{\mathbf{V}}_{f_1\#}^{(n)}(\beta))|\beta \in \mathbb{R}^+\}, \qquad n \in \mathbb{N},$$

with equation

$$\mathring{\text{vech}}(\underset{\sim}{\mathbf{V}}_{f_1\#}^{(n)}(\beta)) := \mathring{\text{vech}}(\mathbf{V}_{\#}^{(n)}) + n^{-1/2}\beta\boldsymbol{\Upsilon}_k(\mathbf{V}_{\#}^{(n)})\underset{\sim}{\boldsymbol{\Delta}}_{f_1}^{(n)}(\mathbf{V}_{\#}^{(n)})$$

(4.2)
$$= \mathring{\text{vech}}(\mathbf{V}_{\#}^{(n)}) + \beta k(k+2)\mathbf{N}_k[\mathbf{I}_{k^2} - (\text{vec }\mathbf{V}_{\#}^{(n)})\mathbf{e}'_{k^2,1}]\text{vec}(\underset{\sim}{\mathbf{W}}_{f_1\#}^{(n)}),$$

where $\mathbf{e}_{k^2,1}$ stands for the first vector of the canonical basis in $\mathbb{R}^{k^2}$ and $\underset{\sim}{\mathbf{W}}_{f_1\#}^{(n)} := \underset{\sim}{\mathbf{W}}_{f_1}^{(n)}(\mathbf{V}_{\#}^{(n)})$; the last equality is obtained exactly as in the proof of Proposition 3.1(iii). Each value of $\beta$ defines on $\mathcal{D}_{\#}^{(n)}$ a sequence of root-$n$ consistent estimators $\underset{\sim}{\mathbf{V}}_{f_1\#}^{(n)}(\beta)$ of $\mathbf{V}$; one of them, namely $\underset{\sim}{\mathbf{V}}_{f_1\#}^{(n)}(\mathcal{J}_k^{-1}(f_1,g_1))$, coincides with $\underset{\sim}{\mathbf{V}}_{f_1\#}^{(n)}$ in (3.1) and is efficient at $\mathcal{P}_{f_1}^{(n)}$ [actually, an estimator $\widehat{\mathbf{V}}^{(n)}$ is efficient iff $\widehat{\mathbf{V}}^{(n)} - \underset{\sim}{\mathbf{V}}_{f_1\#}^{(n)} = o_{\mathrm{P}}(n^{-1/2})$ under $\mathcal{P}_{f_1}^{(n)}$].

However, these estimators $\underset{\sim}{\mathbf{V}}_{f_1\#}^{(n)}(\beta)$ are not locally discrete since the multivariate signs $\mathbf{U}_i^{(n)}$ in $\underset{\sim}{\mathbf{W}}_{f_1\#}^{(n)}$ are not discretized (even though evaluated at $\mathbf{V}_{\#}^{(n)}$); therefore, we discretize them further by discretizing $\underset{\sim}{\mathbf{W}}_{f_1\#}^{(n)}$: let $\underset{\sim}{\mathbf{W}}_{f_1\#\#}^{(n)}$ be the $k \times k$ matrix obtained by mapping each entry $w_{ij\#}^{(n)}$ of $\underset{\sim}{\mathbf{W}}_{f_1\#}^{(n)}$ onto $w_{ij\#\#}^{(n)} := c_1^{-1}\mathrm{sign}(w_{ij\#}^{(n)})n^{-1/2}\lceil n^{1/2}c_1|w_{ij\#}^{(n)}|\rceil$, where $c_1 > 0$ is some arbitrarily large constant. Replacing (4.2) (but keeping the same notation for the sake of simplicity) with

$$\mathrm{v\mathring{e}ch}(\underset{\sim}{\mathbf{V}}_{f_1\#}^{(n)}(\beta_\ell)) := \mathrm{v\mathring{e}ch}(\mathbf{V}_{\#}^{(n)}) + \beta_\ell k(k+2)\mathbf{N}_k[\mathbf{I}_{k^2} - (\mathrm{vec}\,\mathbf{V}_{\#}^{(n)})\mathbf{e}'_{k^2,1}]\mathrm{vec}(\underset{\sim}{\mathbf{W}}_{f_1\#\#}^{(n)})$$

$$(4.3)$$

$$=: \mathrm{v\mathring{e}ch}(\mathbf{V}_{\#}^{(n)}) + n^{-1/2}\beta_\ell\,\boldsymbol{\Upsilon}_k(\mathbf{V}_{\#}^{(n)})\underset{\sim}{\boldsymbol{\Delta}}_{f_1\#}^{(n)}(\mathbf{V}_{\#}^{(n)}), \qquad \ell \in \mathbb{N},$$

where $\beta_\ell := \ell/c_2$, with some other arbitrary constant $c_2 > 0$ yields root-$n$ consistent estimators $\underset{\sim}{\mathbf{V}}_{f_1\#}^{(n)}(\beta_\ell)$ that are locally discrete, in the sense that the number of possible values of $\mathrm{v\mathring{e}ch}(\underset{\sim}{\mathbf{V}}_{f_1\#}^{(n)}(\beta_\ell))$ in balls with $O(n^{-1/2})$ radius centered at $\mathrm{v\mathring{e}ch}(\mathbf{V})$ is bounded as $n \to \infty$. Still for simplicity, we keep the notation $\mathcal{D}_{\#}^{(n)}$ for this new sequence $\mathcal{D}_{\#}^{(n)}(\mathbf{V}_{\#}^{(n)}; \underset{\sim}{\boldsymbol{\Delta}}_{f_1\#}^{(n)}(\mathbf{V}_{\#}^{(n)}))$ of *fully-discretized* half-lines. For any $\ell \in \mathbb{N}$, $\underset{\sim}{\mathbf{V}}_{f_1\#}^{(n)}(\beta_\ell)$ can again serve as the preliminary estimator in a rank-based one-step procedure: letting

$$\mathrm{v\mathring{e}ch}(\underset{\sim}{\mathbf{V}}_{f_1\#}^{(n)}(\beta_\ell;\delta)) := \mathrm{v\mathring{e}ch}(\underset{\sim}{\mathbf{V}}_{f_1\#}^{(n)}(\beta_\ell)) + n^{-1/2}\delta\boldsymbol{\Upsilon}_k(\underset{\sim}{\mathbf{V}}_{f_1\#}^{(n)}(\beta_\ell))\underset{\sim}{\boldsymbol{\Delta}}_{f_1}^{(n)}(\underset{\sim}{\mathbf{V}}_{f_1\#}^{(n)}(\beta_\ell)),$$

$\mathrm{v\mathring{e}ch}(\underset{\sim}{\mathbf{V}}_{f_1\#}^{(n)}(\beta_\ell;\mathcal{J}_k^{-1}(f_1,g_1)))$ is such that

$$(4.4) \qquad \mathrm{v\mathring{e}ch}(\underset{\sim}{\mathbf{V}}_{f_1\#}^{(n)}(\beta_\ell;\mathcal{J}_k^{-1}(f_1,g_1))) - \mathrm{v\mathring{e}ch}(\underset{\sim}{\mathbf{V}}_{f_1\#}^{(n)}) = o_{\mathrm{P}}(n^{-1/2})$$

under $\mathrm{P}_{\sigma^2,\mathbf{V};g_1}^{(n)}$. However, $\mathrm{v\mathring{e}ch}(\underset{\sim}{\mathbf{V}}_{f_1\#}^{(n)}(\beta_\ell;\mathcal{J}_k^{-1}(f_1,g_1)))$ still cannot be computed from the observations.

Denote by $\mathbf{u}_{\mathcal{D}}$ the unit vector along $\mathcal{D}^{(n)}_{\#}$ (corresponding to $\mathcal{D}^{(n)}_{\#}$'s natural orientation as a half-line) and define

$$(4.5) \quad \ell^+ := \min\{\ell \in \mathbb{N}_0 | h^{\#}(\beta_\ell) := \mathbf{u}'_{\mathcal{D}} \mathbf{\Upsilon}(\underset{\sim}{\mathbf{V}}^{(n)}_{f_1\#}(\beta_\ell)) \underset{\sim f_1}{\mathbf{\Delta}^{(n)}}(\underset{\sim}{\mathbf{V}}^{(n)}_{f_1\#}(\beta_\ell)) \leq 0\},$$

$\ell^- := \ell^+ - 1$ and $\beta^\pm := \beta_{\ell^\pm}$. The integers $\ell^\pm$ are random; in order for $\underset{\sim}{\mathbf{V}}^{(n)}_{f_1\#}(\beta^\pm)$ to remain root-$n$ consistent and locally discrete, it is sufficient to check that $\ell^\pm$ is $O_{\mathrm{P}}(1)$. This implies that for any $\varepsilon > 0$, there exist integers $L_\varepsilon$ and $N_\varepsilon$ such that for all $n \geq N_\varepsilon$, the minimization in (4.5) with probability larger than $1 - \varepsilon$ only runs over the finite set $\ell \in \{1,\ldots,L_\varepsilon\}$ (equivalently, over the finite set $\beta \in \{\beta_1,\ldots,\beta_{L_\varepsilon}\}$). In order to show this, let us assume that $\ell^\pm$ is not $O_{\mathrm{P}}(1)$. Then there exists $\varepsilon > 0$ and a sequence $n_i \uparrow \infty$ such that for all $L \in \mathbb{N}$ and some $\sigma^2$, $\underset{\sim}{\mathbf{V}}$ and $g_1$, $\mathrm{P}^{(n_i)}_{\sigma^2,\mathbf{V};g_1}[\ell^- > L] > \varepsilon$.

Pythagoras' Theorem then implies that for $L > c_2 \mathcal{J}_k^{-1}(f_1, g_1)$, with $\mathrm{P}^{(n_i)}_{\sigma^2,\mathbf{V};g_1}$-probability larger than $\varepsilon$,

$$\|\overset{\circ}{\mathrm{vech}}(\underset{\sim}{\mathbf{V}}^{(n_i)}_{f_1\#}(\beta_L; \mathcal{J}_k^{-1}(f_1,g_1))) - \overset{\circ}{\mathrm{vech}}(\underset{\sim}{\mathbf{V}}^{(n_i)}_{f_1\#})\|$$

$$\geq \|\overset{\circ}{\mathrm{vech}}(\underset{\sim}{\mathbf{V}}^{(n_i)}_{f_1\#}(\beta_L)) - \overset{\circ}{\mathrm{vech}}(\underset{\sim}{\mathbf{V}}^{(n_i)}_{f_1\#})\|$$

$$= n_i^{-1/2}(c_2^{-1}L - \mathcal{J}_k^{-1}(f_1,g_1))\|\mathbf{\Upsilon}_k(\mathbf{V}^{(n_i)}_{\#})\underset{\sim f_1\#}{\mathbf{\Delta}^{(n_i)}}(\mathbf{V}^{(n_i)}_{\#})\|$$

which contradicts the fact that (4.4) holds for $\ell = L$. Thus, $\ell^\pm$ are $O_{\mathrm{P}}(1)$ and $\underset{\sim}{\mathbf{V}}^{(n)}_{f_1\#}(\beta^\pm)$ can also serve as initial estimators in a one-step strategy.

The final step in the construction of our estimator $\underset{\sim}{\widehat{\mathbf{V}}}^{(n)}_{f_1\#}$, then, is a "fine tuning" step which consists of selecting an intermediate point between $\beta^-$ and $\beta^+$. This intermediate value, as we shall see, turns out to consistently estimate $\mathcal{J}_k^{-1}(f_1, g_1)$. Denote by $\boldsymbol{\pi}^{(n)}_{\pm}(\delta)$ the projection onto $\mathcal{D}^{(n)}_{\#}$ of $\overset{\circ}{\mathrm{vech}}(\underset{\sim}{\mathbf{V}}^{(n)}_{f_1\#}(\beta^\pm;\delta))$ and let $\pi^{(n)}_{\pm}(\delta) := \|\boldsymbol{\pi}^{(n)}_{\pm}(\delta) - \overset{\circ}{\mathrm{vech}}(\mathbf{V}^{(n)}_{\#})\|$. Note that $\delta \mapsto \pi^{(n)}_{-}(\delta)$ [resp., $\delta \mapsto \pi^{(n)}_{+}(\delta)$] is $\mathcal{P}^{(n)}$-a.e. continuous and strictly monotone increasing (resp., decreasing). Therefore, there exists a unique $\delta^*$ such that $\boldsymbol{\pi}^{(n)}_{-}(\delta^*) = \boldsymbol{\pi}^{(n)}_{+}(\delta^*)$. The proposed $R$-estimator of $\mathbf{V}$ is the shape matrix $\underset{\sim}{\widehat{\mathbf{V}}}^{(n)}_{f_1\#}$ characterized by $\overset{\circ}{\mathrm{vech}}(\underset{\sim}{\widehat{\mathbf{V}}}^{(n)}_{f_1\#}) := \boldsymbol{\pi}^{(n)}_{\pm}(\delta^*)$.

Let us show, to conclude, that $\boldsymbol{\pi}^{(n)}_{\pm}(\delta^*) - \overset{\circ}{\mathrm{vech}}(\underset{\sim}{\mathbf{V}}^{(n)}_{f_1\#}) = o_{\mathrm{P}}(n^{-1/2})$ under $\mathcal{P}^{(n)}_{g_1}$. Either we have $\pi^{(n)}_{-}(\mathcal{J}_k^{-1}(f_1,g_1)) \leq \pi^{(n)}_{+}(\mathcal{J}_k^{-1}(f_1,g_1))$ and

$$\pi^{(n)}_{-}(\mathcal{J}_k^{-1}(f_1,g_1)) \leq \pi^{(n)}_{\pm}(\delta^*) \leq \pi^{(n)}_{+}(\mathcal{J}_k^{-1}(f_1,g_1)),$$

or $\pi_{-}^{(n)}(\mathcal{J}_k^{-1}(f_1, g_1)) > \pi_{+}^{(n)}(\mathcal{J}_k^{-1}(f_1, g_1))$ and

$$\pi_{+}^{(n)}(\mathcal{J}_k^{-1}(f_1, g_1)) < \pi_{\pm}^{(n)}(\delta^*) \leq \pi_{-}^{(n)}(\mathcal{J}_k^{-1}(f_1, g_1)).$$

In both cases, $\boldsymbol{\pi}_{\pm}^{(n)}(\delta^*)$ is in the interval $[\boldsymbol{\pi}_{-}^{(n)}(\mathcal{J}_k^{-1}(f_1, g_1)), \boldsymbol{\pi}_{+}^{(n)}(\mathcal{J}_k^{-1}(f_1, g_1))]$.

Now, both $\boldsymbol{\pi}_{-}^{(n)}(\mathcal{J}_k^{-1}(f_1, g_1))$ and $\boldsymbol{\pi}_{+}^{(n)}(\mathcal{J}_k^{-1}(f_1, g_1))$ are efficient estimators satisfying (3.2) and (3.3). Indeed, from Pythagoras' Theorem,

$$\|\boldsymbol{\pi}_{\pm}^{(n)}(\mathcal{J}_k^{-1}(f_1, g_1)) - \mathring{\text{vech}}(\underset{\sim}{\mathbf{V}}_{f_1\#}^{(n)})\|$$

$$\leq \|\mathring{\text{vech}}(\underset{\sim}{\mathbf{V}}_{f_1\#}^{(n)}(\beta_{\ell\pm}; \mathcal{J}_k^{-1}(f_1, g_1))) - \mathring{\text{vech}}(\underset{\sim}{\mathbf{V}}_{f_1\#}^{(n)})\| = o_{\mathrm{P}}(n^{-1/2})$$

under $\mathcal{P}_{g_1}^{(n)}$. Therefore, as a convex linear combination of $\boldsymbol{\pi}_{-}^{(n)}(\mathcal{J}_k^{-1}(f_1, g_1))$ and $\boldsymbol{\pi}_{+}^{(n)}(\mathcal{J}_k^{-1}(f_1, g_1))$, $\mathring{\text{vech}}(\underset{\sim}{\widehat{\mathbf{V}}}_{f_1\#}^{(n)}) = \boldsymbol{\pi}_{\pm}^{(n)}(\delta^*)$ is also an efficient estimator satisfying (3.2) and (3.3) and, contrary to $\boldsymbol{\pi}_{\pm}^{(n)}(\mathcal{J}_k^{-1}(f_1, g_1))$, it is computable from the sample. Now, clearly,

$$\alpha_\#^* := (\beta_\#^*)^{-1} := [n^{1/2}\|\boldsymbol{\pi}_{\pm}^{(n)}(\delta^*) - \mathring{\text{vech}}(\mathbf{V}_\#^{(n)})\|/\|\boldsymbol{\Upsilon}_k(\mathbf{V}_\#^{(n)})\underset{\sim}{\boldsymbol{\Delta}}_{f_1\#}^{(n)}(\mathbf{V}_\#^{(n)})\|]^{-1}$$

(4.6)

and $(\mathcal{J}_k(f_1)/(\alpha_\#^*)^2)\boldsymbol{\Upsilon}_k(\underset{\sim}{\widehat{\mathbf{V}}}_{f_1\#}^{(n)})$ yield consistent (under $\mathcal{P}_{g_1}^{(n)}$) estimators of $\mathcal{J}_k(f_1, g_1)$ and consistent (under $\mathcal{P}^{(n)}$) estimators the asymptotic covariance matrix of $\mathring{\text{vech}}(\underset{\sim}{\widehat{\mathbf{V}}}_{f_1\#}^{(n)})$, respectively.

4.3. *An original* (*local likelihood*) *method*: *practical implementation*. As usual, the discretization technique which complicates the proofs of asymptotic results and obscures the definition of the estimator makes little sense in practice, where $n$ is fixed. Discretization in the previous sections was achieved in three steps: discretization of Tyler's $\mathbf{V}_T^{(n)}$ into $\mathbf{V}_\#^{(n)}$ (based on $c_0$), discretization of $\underset{\sim}{\boldsymbol{\Delta}}_{f_1}^{(n)}(\mathbf{V}_\#^{(n)})$ into $\underset{\sim}{\boldsymbol{\Delta}}_{f_1\#}^{(n)}(\mathbf{V}_\#^{(n)})$ (based on $c_1$) and discretization of $\beta$ into $\beta_\ell$ (based on $c_2$). The "undiscretized version" $\underset{\sim}{\widehat{\mathbf{V}}}_{f_1}^{(n)}$ of $\underset{\sim}{\widehat{\mathbf{V}}}_{f_1\#}^{(n)}$ corresponds to arbitrarily large values of these three discretization constants, leaving $\mathbf{V}_T^{(n)}$ and $\underset{\sim}{\boldsymbol{\Delta}}_{f_1}^{(n)}$ unchanged and bringing (for the sample size at hand) $\beta^+$ and $\beta^-$ so close to each other that the final tuning [involving the solution $\delta^*$ of $\boldsymbol{\pi}_{-}^{(n)}(\delta) = \boldsymbol{\pi}_{+}^{(n)}(\delta)$] becomes numerically meaningless. Alternatively, denoting by $\underset{\sim}{\widehat{\mathbf{V}}}_{f_1\#}^{(n)}(\mathbf{c})$ the estimator associated with the discretization constants $\mathbf{c} := (c_0, c_1, c_2)$, we have $\widehat{\underset{\sim}{\mathbf{V}}}{}^{(n)}_{f_1} := \lim_{\mathbf{c} \to \infty} \widehat{\underset{\sim}{\mathbf{V}}}{}^{(n)}_{f_1\#}(\mathbf{c})$, where $\mathbf{c} \to \infty$ means that $c_i \to \infty$ for $i = 0, 1, 2$.

This practical implementation $\widehat{\underset{\sim}{\mathbf{V}}}{}^{(n)}_{f_1}$ of $\widehat{\underset{\sim}{\mathbf{V}}}{}^{(n)}_{f_1\#}$ can be obtained more directly as follows. Letting

$$\mathrm{v\mathring{e}ch}(\underset{\sim}{\mathbf{V}}{}^{(n)}_{f_1}(\beta)) := \mathrm{v\mathring{e}ch}(\mathbf{V}^{(n)}_T) + n^{-1/2}\beta \mathbf{\Upsilon}_k(\mathbf{V}^{(n)}_T) \underset{\sim}{\mathbf{\Delta}}{}^{(n)}_{f_1}(\mathbf{V}^{(n)}_T), \qquad \beta \in \mathbb{R}^+,$$

[the undiscretized version of $\mathrm{v\mathring{e}ch}(\underset{\sim}{\mathbf{V}}{}^{(n)}_{f_1\#}(\beta_\ell))$], consider the $\mathcal{P}^{(n)}$-a.e. piecewise continuous function

$$\beta \mapsto h(\beta) := (\underset{\sim}{\mathbf{\Delta}}{}^{(n)}_{f_1}(\mathbf{V}^{(n)}_T))' \mathbf{\Upsilon}_k(\mathbf{V}^{(n)}_T) \mathbf{\Upsilon}_k(\underset{\sim}{\mathbf{V}}{}^{(n)}_{f_1}(\beta)) \underset{\sim}{\mathbf{\Delta}}{}^{(n)}_{f_1}(\underset{\sim}{\mathbf{V}}{}^{(n)}_{f_1}(\beta)),$$

(4.7)
$$\beta \in \mathbb{R}^+,$$

and put $\beta^* := \inf\{\beta > 0 | h(\beta) \leq 0\}$, $\beta^{*-} := \beta^* - 0$ and $\beta^{*+} := \beta^* + 0$. The matrices $\underset{\sim}{\mathbf{V}}{}^{(n)}_{f_1}(\beta^{*-})$ and $\underset{\sim}{\mathbf{V}}{}^{(n)}_{f_1}(\beta^{*+})$ are clearly the "undiscretized counterparts" of $\underset{\sim}{\mathbf{V}}{}^{(n)}_{f_1\#}(\beta^-)$ and $\underset{\sim}{\mathbf{V}}{}^{(n)}_{f_1\#}(\beta^+)$, respectively. However, $\beta \mapsto \underset{\sim}{\mathbf{V}}{}^{(n)}_{f_1}(\beta)$ being continuous, $\underset{\sim}{\mathbf{V}}{}^{(n)}_{f_1}(\beta^{*-}) = \underset{\sim}{\mathbf{V}}{}^{(n)}_{f_1}(\beta^{*+})$. The estimator proposed in Section 4.2 lies between $\underset{\sim}{\mathbf{V}}{}^{(n)}_{f_1\#}(\beta^-)$ and $\underset{\sim}{\mathbf{V}}{}^{(n)}_{f_1\#}(\beta^+)$. Accordingly, the $R$-estimator we are proposing in practice is $\widehat{\underset{\sim}{\mathbf{V}}}{}^{(n)}_{f_1} := \underset{\sim}{\mathbf{V}}{}^{(n)}_{f_1}(\beta^*) = \underset{\sim}{\mathbf{V}}{}^{(n)}_{f_1}(\beta^{*\pm})$; $\alpha^* := (\beta^*)^{-1}$ provides the corresponding estimator of $\mathcal{J}_k(f_1, g_1)$, the "undiscretized" version of (4.6).

Let us stress, however, that all asymptotic properties—including asymptotic optimality—are properties of the discretized estimators $\widehat{\underset{\sim}{\mathbf{V}}}{}^{(n)}_{f_1\#}$, whereas nothing can be said about the asymptotics of the practical implementation $\widehat{\underset{\sim}{\mathbf{V}}}{}^{(n)}_{f_1}$.

**5. Asymptotic affine-equivariance.** An estimator $\mathbf{V}^{(n)}$ of the shape matrix $\mathbf{V}$ is said to be (strictly, that is, for any fixed $n$) affine-equivariant iff for any invertible $k \times k$ matrix $\mathbf{M}$ and any $k$-vector $\mathbf{a}$,

(5.1) $$\mathbf{V}^{(n)}(\mathbf{M}, \mathbf{a}) = (\mathbf{M}\mathbf{V}^{(n)}\mathbf{M}')/(\mathbf{M}\mathbf{V}^{(n)}\mathbf{M}')_{11},$$

where $\mathbf{V}^{(n)}(\mathbf{M}, \mathbf{a})$ denotes the value of the statistic $\mathbf{V}^{(n)}$ computed from the transformed sample $\mathbf{M}\mathbf{X}_1 + \mathbf{a}, \ldots, \mathbf{M}\mathbf{X}_n + \mathbf{a}$. Both Tyler's $\mathbf{V}^{(n)}_T$ and the Gaussian estimator $\mathbf{V}^{(n)}_{\mathcal{G}}$ are affine-equivariant. Unfortunately, the final estimators $\widehat{\underset{\sim}{\mathbf{V}}}{}^{(n)}_{f_1}$ proposed in Section 4.3 are not.

The question arises as to whether $\widehat{\underset{\sim}{\mathbf{V}}}_{f_1}^{(n)}$ is at least *asymptotically* affine-equivariant, that is, whether $\widehat{\underset{\sim}{\mathbf{V}}}_{f_1}^{(n)}$ is asymptotically equivalent to some strictly affine-equivariant sequence (not necessarily a sequence of estimators): for all practical purposes, a sequence of pseudo-estimators, or simply a sequence of random shape matrices, would be fine. Closer inspection of this idea, however, reveals a major conceptual problem. Indeed, recall that all asymptotic results belong to the discretized estimators $\widehat{\underset{\sim}{\mathbf{V}}}_{f_1\#}^{(n)}$, while nothing can be said about the asymptotics of $\widehat{\underset{\sim}{\mathbf{V}}}_{f_1}^{(n)}$: a definition of asymptotic equivariance relying on the asymptotic behavior of $\widehat{\underset{\sim}{\mathbf{V}}}_{f_1}^{(n)}$ is thus totally ineffective.

Therefore, we propose the following, slightly weaker, definition. Denote by

$$\mathcal{S}^{(n)} := \{\mathbf{S}_m^{(n)}(\mathbf{X}^{(n)}) \,|\, m \in \mathbb{N}\} \quad \text{and} \quad \mathcal{T}^{(n)} := \{\mathbf{T}_m^{(n)}(\mathbf{X}^{(n)}) \,|\, m \in \mathbb{N}\}, \qquad n \in \mathbb{N},$$

two countable sequences of $\mathbf{X}^{(n)}$-measurable random vectors or matrices such that the a.s. limits $\mathbf{S}^{(n)} := \lim_{m \to \infty} \mathbf{S}_m^{(n)}(\mathbf{X}^{(n)})$ and $\mathbf{T}^{(n)} := \lim_{m \to \infty} \mathbf{T}_m^{(n)}(\mathbf{X}^{(n)})$ exist for all fixed $n$. Then if

(i) $\mathcal{S}^{(n)}$ and $\mathcal{T}^{(n)}$ are asymptotically equivalent, meaning that for all $m$ (or, more generally, for $m$ large enough), $\mathbf{S}_m^{(n)}(\mathbf{X}^{(n)}) - \mathbf{T}_m^{(n)}(\mathbf{X}^{(n)}) = o_{\mathrm{P}}(n^{-1/2})$ as $n \to \infty$, and if
(ii) $\mathbf{S}^{(n)}$ is strictly equivariant,

we may consider that $\mathbf{T}^{(n)}$ inherits, under approximate or asymptotic form, the equivariance property of $\mathbf{S}^{(n)}$; we say that $\mathbf{T}^{(n)}$ is *weakly asymptotically equivariant*.

In order to show that the proposed estimators $\widehat{\underset{\sim}{\mathbf{V}}}_{f_1}^{(n)} := \lim_{\mathbf{c} \to \infty} \widehat{\underset{\sim}{\mathbf{V}}}_{f_1\#}^{(n)}(\mathbf{c})$ are weakly asymptotically affine-equivariant, consider the class $\mathcal{T}^{(n)} := \{\widehat{\underset{\sim}{\mathbf{V}}}_{f_1\#}^{(n)}(\mathbf{c}_m) \,|\, m \in \mathbb{N}\}$, where the sequence $\mathbf{c}_m = (c_{m,0}, c_{m,1}, c_{m,2})$ is such that $\lim_{m \to \infty} c_{m,i} = \infty$, $i = 0, 1, 2$ and let us construct a class $\mathcal{S}^{(n)}$ such that conditions (i) and (ii) for weak asymptotic equivariance are satisfied. Incidentally, note that a choice of the form $\mathcal{S}^{(n)} := \{\underset{\sim}{\mathbf{V}}_{f_1\#}^{(n)}(c_{0,m}) \,|\, m \in \mathbb{N}\}$ (with $c_{0,m} \to \infty$), where $\underset{\sim}{\mathbf{V}}_{f_1\#}^{(n)}(c_0)$ denotes the pseudo-estimator defined in (3.1), is not suitable since the corresponding practical implementation $\underset{\sim}{\mathbf{V}}_f^{(n)} := \lim_{c_0 \to \infty} \underset{\sim}{\mathbf{V}}_{f_1\#}^{(n)}(c_0)$ is not strictly affine-equivariant.

Inspired by $\widehat{\underset{\sim}{\mathbf{V}}}_{f_1\#}^{(n)}$'s representation (3.9) as a linear combination of $\mathbf{V}_{\#}^{(n)}$ and the rank-based shape matrix $\underset{\sim}{\mathbf{W}}_{f_1\#}^{(n)}$ defined in (3.10), consider now the

shape pseudo-estimators

(5.2) $$\underset{\approx}{\mathbf{V}}_{f_1\#}^{(n)} = \underset{\approx}{\mathbf{V}}_{f_1\#}^{(n)}(c_0) := \underset{\sim}{\mathbf{B}}_{f_1\#}^{(n)} / (\underset{\sim}{\mathbf{B}}_{f_1\#}^{(n)})_{11}$$

with $\underset{\sim}{\mathbf{B}}_{f_1\#}^{(n)} := (1 - \frac{k(k+2)}{\mathcal{J}_k(f_1,g_1)}) \mathbf{V}_\#^{(n)} + \frac{k(k+2)}{\mathcal{J}_k(f_1,g_1)} \underset{\sim}{\mathbf{W}}_{f_1\#}^{(n)}$, where $c_0$ is the constant used in the discretization of Tyler's $\mathbf{V}_T^{(n)}$. Although, due to discretization, neither $\mathbf{V}_\#^{(n)}$ nor $\underset{\approx}{\mathbf{V}}_{f_1\#}^{(n)}$ are affine-equivariant for fixed $n$, the class $\mathcal{S}^{(n)} := \{\underset{\approx}{\mathbf{V}}_{f_1\#}^{(n)}(c_{0,m}) \mid m \in \mathbb{N}\}$ allows us to establish the weak asymptotic affine-equivariance of $\widehat{\underset{\sim}{\mathbf{V}}}_{f_1}^{(n)}$, as shown in the following proposition.

PROPOSITION 5.1. *Denote by* $\underset{\approx}{\mathbf{V}}_{f_1\#}^{(n)} := \underset{\approx}{\mathbf{V}}_{f_1\#}^{(n)}(c_0)$ *and by* $\widehat{\underset{\sim}{\mathbf{V}}}_{f_1\#}^{(n)} := \widehat{\underset{\sim}{\mathbf{V}}}_{f_1\#}^{(n)}(\mathbf{c})$ *the pseudo-estimator defined in* (5.2) *and the estimator defined in Section* 4.2, *respectively. Then* (i) $\underset{\approx}{\mathbf{V}}_{f_1\#}^{(n)} - \widehat{\underset{\sim}{\mathbf{V}}}_{f_1\#}^{(n)} = o_P(n^{-1/2})$ *under* $\mathcal{P}^{(n)}$ *as* $n \to \infty$ *and* (ii) *the practical implementation* $\underset{\approx}{\mathbf{V}}_{f_1}^{(n)} := \lim_{m\to\infty} \underset{\approx}{\mathbf{V}}_{f_1\#}^{(n)}(c_{0,m})$ *is strictly affine-equivariant.*

PROOF. See Section A.3. □

Whether or not weak asymptotic equivariance is a satisfactory property is a matter of statistical taste. If it is, then this section shows that $\widehat{\underset{\sim}{\mathbf{V}}}_{f_1}^{(n)}$ is the estimator to be used. The reader who feels that strict equivariance is an essential requirement is referred to [15], where it is shown that an adequate modification of $\widehat{\underset{\sim}{\mathbf{V}}}_{f_1}^{(n)}$ producing a strictly equivariant $\widehat{\underset{\approx}{\mathbf{V}}}_{f_1}^{(n)}$ is possible (at the price of some technicalities). Alternatively, it is shown in [8] that, under mild additional assumptions, an affine-equivariant $R$-estimator of shape also can be obtained from iterating the mapping $\mathbf{V} \mapsto \underset{\sim}{\mathbf{W}}_{f_1}^{(n)}(\mathbf{V}) / (\underset{\sim}{\mathbf{W}}_{f_1}^{(n)}(\mathbf{V}))_{11}$, where $\underset{\sim}{\mathbf{W}}_{f_1}^{(n)}$ is defined in (3.10).

**6. Simulations.** In this section, we conduct a Monte Carlo study in order to compare the finite-sample performances of the one-step $R$-estimators $\widehat{\underset{\sim}{\mathbf{V}}}_{f_1}^{(n)}$ proposed in Section 4.3 (as well as those of their analogs using the Gaussian estimator $\mathbf{V}_\mathcal{G}^{(n)}$, instead of Tyler's $\mathbf{V}_T^{(n)}$, as a preliminary estimator) to those of $\mathbf{V}_T^{(n)}$ and $\mathbf{V}_\mathcal{G}^{(n)}$ themselves. We restrict our attention to the bivariate spherical case ($\mathbf{V} = \mathbf{I}_2$). We generated $M = 1{,}000$ samples of i.i.d. observations $\mathbf{X}_1, \ldots, \mathbf{X}_n$, with sizes $n = 50$ and $n = 250$, from the bivariate standard normal distribution ($\mathcal{N}$), the Student distributions ($t_{0.5}$), ($t_3$) and ($t_{10}$) (with 0.5, 3 and 10 degrees of freedom) and the power-exponential distributions ($e_3$) and ($e_5$) (with parameters $\eta = 3$ and 5); for details on power-exponential densities, see Section HP1.2. This choice of Student and power-exponential distributions allows for the consideration of heavier-than-normal and lighter-than-normal tail distributions, respectively.

TABLE 2
*Empirical bias and mean-square error, under various bivariate t-, power-exponential and normal densities, of the preliminary estimators $\mathbf{V}_{\mathcal{G}}^{(n)}$ and $\mathbf{V}_{T}^{(n)}$ and the corresponding one-step R-estimators $\widehat{\mathbf{V}}_{0.5}^{(n)}$, $\widehat{\mathbf{V}}_{3}^{(n)}$, $\widehat{\mathbf{V}}_{10}^{(n)}$ and $\widehat{\mathbf{V}}_{\mathrm{vdW}}^{(n)}$. The simulation is based on 1000 replications; sample size is $n=50/n=250$.*

| | Preliminary estimator | BIAS ($n=50/n=250$) | | | | | |
|---|---|---|---|---|---|---|---|
| | | $t_{0.5}$ | $t_3$ | $t_{10}$ | $\mathcal{N}$ | $e_3$ | $e_5$ |
| – | $\mathbf{V}_{T}^{(n)}$ | 0.0042/−0.0043 | −0.0038/−0.0043 | −0.0016/0.0003 | 0.0006/−0.0030 | 0.0067/0.0070 | −0.0070/−0.0023 |
| | | 0.0830/0.0207 | 0.0973/0.0219 | 0.0865/0.0062 | 0.0895/0.0024 | 0.1118/0.0201 | 0.0906/0.0072 |
| – | $\mathbf{V}_{\mathcal{G}}^{(n)}$ | −0.6148/−0.0522 | 0.0012/−0.0005 | −0.0003/−0.0010 | −0.0058/0.0005 | 0.0025/0.0021 | −0.0024/−0.0021 |
| | | 310.8334/20.6781 | 0.1782/0.0410 | 0.0497/0.0058 | 0.0375/0.0024 | 0.0484/0.0041 | 0.0308/0.0006 |
| $\widehat{\mathbf{V}}_{0.5}^{(n)}$ | $\mathbf{V}_{T}^{(n)}$ | 0.0034/−0.0024 | −0.0004/−0.0031 | 0.0004/−0.0006 | −0.0006/−0.0019 | 0.0039/0.0043 | −0.0030/−0.0026 |
| | | 0.0771/0.0183 | 0.0806/0.0180 | 0.0619/0.0031 | 0.0674/0.0030 | 0.0821/0.0115 | 0.0664/0.0037 |
| | $\mathbf{V}_{\mathcal{G}}^{(n)}$ | 0.0001/0.0021 | 0.0004/−0.0030 | −0.0005/−0.0006 | −0.0007/−0.0019 | 0.0033/0.0043 | −0.0036/−0.0026 |
| | | 0.0798/0.0171 | 0.0782/0.0178 | 0.0612/0.0032 | 0.0671/0.0032 | 0.0820/0.0116 | 0.0661/0.0037 |
| $\widehat{\mathbf{V}}_{3}^{(n)}$ | $\mathbf{V}_{T}^{(n)}$ | 0.0002/−0.0014 | 0.0019/−0.0017 | 0.0005/−0.0009 | −0.0024/−0.0004 | 0.0023/0.0022 | −0.0017/−0.0022 |
| | | 0.0861/0.0216 | 0.0680/0.0142 | 0.0438/0.0024 | 0.0444/0.0028 | 0.0533/0.0047 | 0.0338/0.0006 |
| | $\mathbf{V}_{\mathcal{G}}^{(n)}$ | 0.0014/0.0051 | 0.0028/−0.0017 | 0.0002/−0.0009 | −0.0021/−0.0004 | 0.0023/0.0023 | −0.0019/−0.0021 |
| | | 0.1717/0.0219 | 0.0665/0.0140 | 0.0433/0.0023 | 0.0442/0.0030 | 0.0531/0.0043 | 0.0336/0.0006 |
| $\widehat{\mathbf{V}}_{10}^{(n)}$ | $\mathbf{V}_{T}^{(n)}$ | −0.0001/−0.0008 | 0.0025/−0.0015 | 0.0004/−0.0008 | −0.0036/0.0001 | 0.0023/0.0014 | −0.0019/−0.0021 |
| | | 0.0962/0.0250 | 0.0681/−0.0261 | 0.0427/0.0029 | 0.0395/0.0026 | 0.0441/0.0032 | 0.0253/0.0000 |
| | $\mathbf{V}_{\mathcal{G}}^{(n)}$ | 0.0037/0.0075 | 0.0034/−0.0014 | 0.0001/−0.0008 | −0.0031/0.0001 | 0.0023/0.0016 | −0.0019/−0.0020 |
| | | 0.1074/0.0254 | 0.0672/0.0128 | 0.0419/0.0028 | 0.0398/0.0028 | 0.0440/0.0032 | 0.0250/−0.0000 |
| $\widehat{\mathbf{V}}_{\mathrm{vdW}}^{(n)}$ | $\mathbf{V}_{T}^{(n)}$ | 0.0005/−0.0003 | 0.0027/−0.0014 | 0.0005/−0.0007 | −0.0044/0.0003 | 0.0024/0.0011 | −0.0024/−0.0020 |
| | | 0.1057/0.0281 | 0.0702/0.0124 | 0.0441/0.0036 | 0.0387/0.0025 | 0.0404/0.0026 | 0.0217/−0.0000 |
| | $\mathbf{V}_{\mathcal{G}}^{(n)}$ | 0.0034/0.0091 | 0.0035/−0.0013 | −0.0001/−0.0007 | −0.0041/0.0004 | 0.0024/0.0013 | −0.0022/−0.0019 |
| | | 0.1164/0.0284 | 0.0696/0.0122 | 0.0435/0.0036 | 0.0392/0.0026 | 0.0402/0.0026 | 0.0211/−0.0001 |

TABLE 2
*(Continued).*

|  |  | MSE ($n=50/n=250$) | | | | | |
|---|---|---|---|---|---|---|---|
|  |  | $t_{0.5}$ | $t_3$ | $t_{10}$ | $\mathcal{N}$ | $e_3$ | $e_5$ |
| – | $\mathbf{V}_T^{(n)}$ | 0.0410/0.0083 | 0.0407/0.0081 | 0.0408/0.0075 | 0.0404/0.0075 | 0.0444/0.0080 | 0.0423/0.0085 |
|  |  | 0.2009/0.0392 | 0.2467/0.0357 | 0.2192/0.0337 | 0.2311/0.0369 | 0.2163/0.0337 | 0.2031/0.0320 |
| – | $\mathbf{V}_\mathcal{G}^{(n)}$ | 298.8463/11.3416 | 0.1033/0.0329 | 0.0265/0.0050 | 0.0183/0.0038 | 0.0155/0.0028 | 0.0138/0.0029 |
|  |  | 80,313,350/42,948 | 0.7141/0.2358 | 0.1247/0.0211 | 0.0941/0.0175 | 0.0624/0.0115 | 0.0617/0.0109 |
| $\widehat{\mathbf{V}}_{0.5}^{(n)}$ | $\mathbf{V}_T^{(n)}$ | 0.0368/0.0075 | 0.0328/0.0065 | 0.0312/0.0058 | 0.0307/0.0058 | 0.0320/0.0057 | 0.0296/0.0061 |
|  |  | 0.1862/0.0339 | 0.1879/0.0285 | 0.1629/0.0258 | 0.1701/0.0282 | 0.1425/0.0233 | 0.1411/0.0223 |
|  | $\mathbf{V}_\mathcal{G}^{(n)}$ | 0.1152/0.0278 | 0.0337/0.0065 | 0.0308/0.0057 | 0.0309/0.0058 | 0.0318/0.0057 | 0.0294/0.0061 |
|  |  | 0.2700/0.0566 | 0.1852/0.0284 | 0.1614/0.0258 | 0.1686/0.0281 | 0.1416/0.0233 | 0.1398/0.0223 |
| $\widehat{\mathbf{V}}_3^{(n)}$ | $\mathbf{V}_T^{(n)}$ | 0.0419/0.0090 | 0.0290/0.0057 | 0.0238/0.0044 | 0.0208/0.0042 | 0.0178/0.0031 | 0.0149/0.0030 |
|  |  | 0.2239/0.0371 | 0.1546/0.0247 | 0.1169/0.0199 | 0.1138/0.0198 | 0.0715/0.0127 | 0.0676/0.0112 |
|  | $\mathbf{V}_\mathcal{G}^{(n)}$ | 0.1184/0.0295 | 0.0296/0.0058 | 0.0235/0.0044 | 0.0209/0.0042 | 0.0175/0.0031 | 0.0146/0.0030 |
|  |  | 5.6092/0.0598 | 0.1537/0.0247 | 0.1162/0.0199 | 0.1132/0.0197 | 0.0709/0.0128 | 0.0668/0.0112 |
| $\widehat{\mathbf{V}}_{10}^{(n)}$ | $\mathbf{V}_T^{(n)}$ | 0.0490/0.0106 | 0.0300/0.0060 | 0.0234/0.0043 | 0.0191/0.0039 | 0.0147/0.0025 | 0.0118/0.0022 |
|  |  | 0.2701/0.0428 | 0.1579/1.5539 | 0.1117/0.0191 | 0.1005/0.0180 | 0.0568/0.0102 | 0.0519/0.0084 |
|  | $\mathbf{V}_\mathcal{G}^{(n)}$ | 0.1307/0.0339 | 0.0306/0.0060 | 0.0232/0.0043 | 0.0190/0.0039 | 0.0143/0.0025 | 0.0114/0.0022 |
|  |  | 0.3796/0.0662 | 0.1583/0.0253 | 0.1108/0.0191 | 0.1006/0.0180 | 0.0562/0.0101 | 0.0511/0.0083 |
| $\widehat{\mathbf{V}}_{\text{vdW}}^{(n)}$ | $\mathbf{V}_T^{(n)}$ | 0.0552/0.0121 | 0.0316/0.0064 | 0.0238/0.0044 | 0.0187/0.0039 | 0.0135/0.0022 | 0.0106/0.0019 |
|  |  | 0.3134/0.0486 | 0.1652/0.0267 | 0.1129/0.0192 | 0.0964/0.0176 | 0.0518/0.0092 | 0.0457/0.0073 |
|  | $\mathbf{V}_\mathcal{G}^{(n)}$ | 0.1406/0.0377 | 0.0322/0.0064 | 0.0238/0.0044 | 0.0185/0.0039 | 0.0131/0.0022 | 0.0102/0.0018 |
|  |  | 0.4237/0.0726 | 0.1665/0.0266 | 0.1121/0.0192 | 0.0967/0.0175 | 0.0511/0.0092 | 0.0449/0.0072 |

For each replication, we computed $\mathbf{V}_T^{(n)}$, $\mathbf{V}_\mathcal{G}^{(n)}$ and the $\mathbf{V}_T^{(n)}$- and $\mathbf{V}_\mathcal{G}^{(n)}$-based one-step $R$-estimators $\widehat{\underset{\sim}{\mathbf{V}}}_{\mathrm{vdW}}^{(n)}$, $\widehat{\underset{\sim}{\mathbf{V}}}_{0.5}^{(n)}$, $\widehat{\underset{\sim}{\mathbf{V}}}_{3}^{(n)}$ and $\widehat{\underset{\sim}{\mathbf{V}}}_{10}^{(n)}$, corresponding to semiparametric efficiency at Gaussian and Student densities with 0.5, 3 and 10 degrees of freedom, respectively. In Table 2, we report, for each estimate $\mathbf{V}^{(n)}(l) = (V_{ij}^{(n)}(l))$, the two components of the average bias

$$\mathrm{BIAS}^{(n)} := \frac{1}{M} \sum_{l=1}^{M} \overset{\circ}{\mathrm{vech}}(\mathbf{V}^{(n)}(l) - \mathbf{V}) = \frac{1}{M} \sum_{l=1}^{M} (V_{12}^{(n)}(l), V_{22}^{(n)}(l) - 1)'$$

and the two components of the mean square error

$$\mathrm{MSE}^{(n)} := \frac{1}{M} \sum_{i=1}^{M} ((V_{12}^{(n)}(l))^2, (V_{22}^{(n)}(l) - 1)^2)',$$

based on the $M$ replications $\mathbf{V}^{(n)}(l)$, $l = 1, \ldots, M$.

These simulations show that the proposed rank-based estimators behave remarkably well under all distributions under consideration and significantly improve on Tyler's estimator. They confirm the optimality of the Tyler-based $f_1$-score $R$-estimators under radial density $f$ and essentially agree with the ARE rankings presented in Table 1. Also, the van der Waerden rank-based estimator (based on preliminary estimator $\mathbf{V}_T^{(n)}$ or $\mathbf{V}_\mathcal{G}^{(n)}$) uniformly dominates the parametric Gaussian estimator $\mathbf{V}_\mathcal{G}^{(n)}$ and performs equally well in the normal case; this dominance, which is observed under both lighter-than-normal and heavier-than-normal tail distributions, provides an empirical validation of the Chernoff–Savage result of [38].

The behavior of one-step rank-based estimators does not seem to depend much on the preliminary estimator used ($\mathbf{V}_T^{(n)}$ or $\mathbf{V}_\mathcal{G}^{(n)}$), confirming that the influence of the preliminary estimator is asymptotically nil. More surprising is the fact that $R$-estimators based on $\mathbf{V}_\mathcal{G}^{(n)}$ behave reasonably well under heavy tails (under $t_{0.5}$), although $\mathbf{V}_\mathcal{G}^{(n)}$ itself is not even root-$n$ consistent there (which explains its total collapse under $t_{0.5}$). Quite remarkably, these conclusions are equally valid for small ($n = 50$) as for large ($n = 250$) sample sizes. This is another non-negligible advantage of our method over kernel-based ones (see Section 4.1), which typically require much larger sample sizes.

## APPENDIX

**A.1. Local asymptotic linearity.** Rather than Proposition 2.1(v), we prove in this section a more general asymptotic linearity result in which both the location and the shape parameters are locally perturbed.

PROPOSITION A.1. *For any bounded sequences of $k$-dimensional vectors $\mathbf{t}^{(n)}$ and symmetric matrices $\mathbf{v}^{(n)}$ satisfying $v_{11}^{(n)} = 0$ and for any $g_1 \in \mathcal{F}_A$, the central sequence $\underset{\sim}{\boldsymbol{\Delta}}_{f_1}^{(n)}(\boldsymbol{\theta}, \mathbf{V})$ satisfies, under $\mathrm{P}_{\boldsymbol{\theta}, \sigma^2, \mathbf{V}; g_1}^{(n)}$, as $n \to \infty$, the asymptotic linearity property*

$$\begin{aligned}
\underset{\sim}{\boldsymbol{\Delta}}_{f_1}^{(n)}(\boldsymbol{\theta} + n^{-1/2}\mathbf{t}^{(n)}, \mathbf{V} + n^{-1/2}\mathbf{v}^{(n)}) &- \underset{\sim}{\boldsymbol{\Delta}}_{f_1}^{(n)}(\boldsymbol{\theta}, \mathbf{V}) \\
&= -\boldsymbol{\Gamma}_{f_1, g_1}^*(\mathbf{V})\overset{\circ}{\mathrm{vech}}(\mathbf{v}^{(n)}) + o_\mathrm{P}(1).
\end{aligned}$$
(A.1)

The proof of Proposition A.1 relies on a series of lemmas. In this section, we let $\boldsymbol{\theta}^n := \boldsymbol{\theta} + n^{-1/2}\mathbf{t}^{(n)}$ and $\mathbf{V}^n := \mathbf{V} + n^{-1/2}\mathbf{v}^{(n)}$. Accordingly, let $\mathbf{Z}_i^0 := \mathbf{V}^{-1/2}(\mathbf{X}_i - \boldsymbol{\theta})$, $d_i^0 := \|\mathbf{Z}_i^0\|$, $\mathbf{U}_i^0 := \mathbf{Z}_i^0/d_i^0$, $\mathbf{Z}_i^n := (\mathbf{V}^n)^{-1/2}(\mathbf{X}_i - \boldsymbol{\theta}^n)$, $d_i^n := \|\mathbf{Z}_i^n\|$ and $\mathbf{U}_i^n := \mathbf{Z}_i^n/d_i^n$. We begin with the following preliminary result:

LEMMA A.1. *For all $i$, as $n \to \infty$, under $\mathrm{P}_{\boldsymbol{\theta}, \sigma^2, \mathbf{V}; g_1}^{(n)}$,*

(i) $|d_i^n - d_i^0| = o_\mathrm{P}(1)$ *and*
(ii) $\|\mathbf{U}_i^n - \mathbf{U}_i^0\| = o_\mathrm{P}(1)$.

PROOF. First, note that, defining $\|\mathbf{M}\|_\mathcal{L} := \sup_{\|\mathbf{x}\|=1} \|\mathbf{M}\mathbf{x}\|$,

$$\begin{aligned}
\|\mathbf{Z}_i^n - \mathbf{Z}_i^0\| &\leq \|(\mathbf{V}^n)^{-1/2}(\boldsymbol{\theta} - \boldsymbol{\theta}^n)\| + \|((\mathbf{V}^n)^{-1/2} - \mathbf{V}^{-1/2})(\mathbf{X}_i - \boldsymbol{\theta})\| \\
&\leq n^{-1/2}\|(\mathbf{V}^n)^{-1/2}\|_\mathcal{L}\|\mathbf{t}^{(n)}\| + \|(\mathbf{V}^n)^{-1/2} - \mathbf{V}^{-1/2}\|_\mathcal{L}\|\mathbf{V}^{1/2}\|_\mathcal{L} d_i^0 \\
&\leq C(n)(1 + d_i^0)
\end{aligned}$$

for some positive sequence $C(n)$, with $C(n) = o(1)$ as $n \to \infty$. Now, since for all $\delta > 0$, $\mathrm{P}_{\boldsymbol{\theta}, \sigma^2, \mathbf{V}; g_1}^{(n)}[C(n)(d_i^0)^a > \delta] = o(1)$ as $n \to \infty$ ($a = -1, 0, 1$), we obtain that $\|\mathbf{Z}_i^n - \mathbf{Z}_i^0\|$ and $\|\mathbf{Z}_i^n - \mathbf{Z}_i^0\|/d_i^0$ are $o_\mathrm{P}(1)$ under $\mathrm{P}_{\boldsymbol{\theta}, \sigma^2, \mathbf{V}; g_1}^{(n)}$ as $n \to \infty$. The result follows since (i) $|d_i^n - d_i^0| \leq \|\mathbf{Z}_i^n - \mathbf{Z}_i^0\|$ and (ii) $\|\mathbf{U}_i^n - \mathbf{U}_i^0\| \leq |(1/d_i^n - 1/d_i^0)|\|\mathbf{Z}_i^n\| + \|\mathbf{Z}_i^n - \mathbf{Z}_i^0\|/d_i^0 \leq 2\|\mathbf{Z}_i^n - \mathbf{Z}_i^0\|/d_i^0$. □

PROOF OF PROPOSITION A.1. We first consider the following truncation of the score function $K_{f_1}$. For all $\ell \in \mathbb{N}_0$, define

$$K_{f_1}^{(\ell)}(u) := K_{f_1}\left(\frac{2}{\ell}\right)\ell\left(u - \frac{1}{\ell}\right)I_{[\frac{1}{\ell}<u\leq\frac{2}{\ell}]} + K_{f_1}(u)I_{[\frac{2}{\ell}<u\leq 1-\frac{2}{\ell}]}$$
$$+ K_{f_1}\left(1 - \frac{2}{\ell}\right)\ell\left(\left(1 - \frac{1}{\ell}\right) - u\right)I_{[1-\frac{2}{\ell}<u\leq 1-\frac{1}{\ell}]},$$

where $I_A$ denotes the indicator function of $A$. Since $u \mapsto K_{f_1}(u)$ is continuous, the functions $u \mapsto K_{f_1}^{(\ell)}(u)$ are also continuous on $(0,1)$. It follows that the truncated scores $K_{f_1}^{(\ell)}$ are bounded for all $\ell$. Clearly, it can be safely assumed that $K_{f_1}$ is a monotone increasing function (rather than the difference of two monotone increasing functions) so that there exists some $L$ such that $|K_{f_1}^{(\ell)}(u)| \leq |K_{f_1}(u)|$ for all $u \in (0,1)$ and all $\ell \geq L$.

We have to prove that, under $P_{\boldsymbol{\theta},\sigma^2,\mathbf{V};g_1}^{(n)}$, as $n \to \infty$,

(A.2) $\underset{\sim}{\boldsymbol{\Delta}}_{f_1}^{(n)}(\boldsymbol{\theta}^n, \mathbf{V}^n) - \underset{\sim}{\boldsymbol{\Delta}}_{f_1}^{(n)}(\boldsymbol{\theta}, \mathbf{V}) + \mathcal{J}_k(f_1,g_1)\boldsymbol{\Upsilon}_k^{-1}(\mathbf{V})\overset{\circ}{\text{vech}}(\mathbf{v}^{(n)})$

is $o_P(1)$. Proposition 2.1(ii) shows that $\underset{\sim}{\boldsymbol{\Delta}}_{f_1}^{(n)}(\boldsymbol{\theta}, \mathbf{V}) - \boldsymbol{\Delta}_{f_1,g_1}^{*(n)}(\boldsymbol{\theta}, \mathbf{V})$ is $o_P(1)$ as $n \to \infty$, under the same sequence of hypotheses. Similarly, the difference $\underset{\sim}{\boldsymbol{\Delta}}_{f_1}^{(n)}(\boldsymbol{\theta}^n, \mathbf{V}^n) - \boldsymbol{\Delta}_{f_1,g_1}^{*(n)}(\boldsymbol{\theta}^n, \mathbf{V}^n)$ is $o_P(1)$ as $n \to \infty$, under $P_{\boldsymbol{\theta}^n,\sigma^2,\mathbf{V}^n;g_1}^{(n)}$ (hence, from contiguity, also under $P_{\boldsymbol{\theta},\sigma^2,\mathbf{V};g_1}^{(n)}$). Consequently, (A.2) is asymptotically equivalent, under $P_{\boldsymbol{\theta},\sigma^2,\mathbf{V};g_1}^{(n)}$, to

(A.3) $\boldsymbol{\Delta}_{f_1,g_1}^{*(n)}(\boldsymbol{\theta}^n, \mathbf{V}^n) - \boldsymbol{\Delta}_{f_1,g_1}^{*(n)}(\boldsymbol{\theta}, \mathbf{V}) + \mathcal{J}_k(f_1,g_1)\boldsymbol{\Upsilon}_k^{-1}(\mathbf{V})\overset{\circ}{\text{vech}}(\mathbf{v}^{(n)}).$

Now, $n^{-1/2}\mathbf{J}_k^\perp \text{vec}\left[\sum_{i=1}^n K_{f_1}(\tilde{G}_{1k}(d_i^n/\sigma))\mathbf{U}_i^n\mathbf{U}_i^{n\prime}\right]$, under $P_{\boldsymbol{\theta}^n,\sigma^2,\mathbf{V}^n;g_1}^{(n)}$, is asymptotically normal as $n \to \infty$, with mean zero and covariance matrix $(k(k+2))^{-1}\mathcal{J}_k(f_1)[\mathbf{I}_{k^2} + \mathbf{K}_k - \frac{2}{k}\mathbf{J}_k]$ so that

$$\frac{1}{2}n^{-1/2}\mathbf{M}_k\left[((\mathbf{V}^n)^{\otimes 2})^{-1/2} - (\mathbf{V}^{\otimes 2})^{-1/2}\right]\mathbf{J}_k^\perp \text{vec}\left[\sum_{i=1}^n K_{f_1}(\tilde{G}_{1k}(d_i^n/\sigma))\mathbf{U}_i^n\mathbf{U}_i^{n\prime}\right]$$

is $o_P(1)$ as $n \to \infty$ under $P_{\boldsymbol{\theta}^n,\sigma^2,\mathbf{V}^n;g_1}^{(n)}$, as well as under $P_{\boldsymbol{\theta},\sigma^2,\mathbf{V};g_1}^{(n)}$ (by contiguity). Consequently, (A.3) is asymptotically equivalent, under $P_{\boldsymbol{\theta},\sigma^2,\mathbf{V};g_1}^{(n)}$, to

$$\mathbf{C}^{(n)} := \frac{1}{2}n^{-1/2}\mathbf{M}_k(\mathbf{V}^{\otimes 2})^{-1/2}\mathbf{J}_k^\perp \text{vec}\left[\sum_{i=1}^n K_{f_1}(\tilde{G}_{1k}(d_i^n/\sigma))\mathbf{U}_i^n\mathbf{U}_i^{n\prime}\right]$$

(A.4) $\quad -\frac{1}{2}n^{-1/2}\mathbf{M}_k(\mathbf{V}^{\otimes 2})^{-1/2}\mathbf{J}_k^\perp \text{vec}\left[\sum_{i=1}^n K_{f_1}(\tilde{G}_{1k}(d_i^0/\sigma))\mathbf{U}_i^0\mathbf{U}_i^{0\prime}\right]$

$$+ \mathcal{J}_k(f_1, g_1) \Upsilon_k^{-1}(\mathbf{V}) \mathring{\text{vech}}(\mathbf{v}^{(n)})$$

and we need only prove that $\mathbf{C}^{(n)} = o_\mathrm{P}(1)$. Decompose $\mathbf{C}^{(n)}$ into $\mathbf{C}^{(n)} = \mathbf{D}_1^{(n;\ell)} + \mathbf{D}_2^{(n;\ell)} - \mathbf{R}_1^{(n;\ell)} + \mathbf{R}_2^{(n;\ell)} + \mathbf{R}_3^{(n;\ell)}$, where, denoting by $\mathrm{E}_0$ expectation under $\mathrm{P}_{\boldsymbol{\theta},\sigma^2,\mathbf{V};g_1}^{(n)}$ and defining $\mathcal{J}_k^{(\ell)}(f_1; g_1) := \int_0^1 K_{f_1}^{(\ell)}(u) K_{g_1}(u)\, du$,

$$\mathbf{D}_1^{(n;\ell)} := \frac{1}{2} n^{-1/2} \mathbf{M}_k(\mathbf{V}^{\otimes 2})^{-1/2} \mathbf{J}_k^\perp \text{vec}\left[\sum_{i=1}^n K_{f_1}^{(\ell)}(\tilde{G}_{1k}(d_i^n/\sigma)) \mathbf{U}_i^n \mathbf{U}_i^{n\prime}\right]$$

$$- \frac{1}{2} n^{-1/2} \mathbf{M}_k(\mathbf{V}^{\otimes 2})^{-1/2} \mathbf{J}_k^\perp \text{vec}\left[\sum_{i=1}^n K_{f_1}^{(\ell)}(\tilde{G}_{1k}(d_i^0/\sigma)) \mathbf{U}_i^0 \mathbf{U}_i^{0\prime}\right]$$

$$- \frac{1}{2} n^{-1/2} \mathbf{M}_k(\mathbf{V}^{\otimes 2})^{-1/2} \mathbf{J}_k^\perp \mathrm{E}_0\left[\text{vec}\left[\sum_{i=1}^n K_{f_1}^{(\ell)}(\tilde{G}_{1k}(d_i^n/\sigma)) \mathbf{U}_i^n \mathbf{U}_i^{n\prime}\right]\right],$$

$$\mathbf{D}_2^{(n;\ell)} := \frac{1}{2} n^{-1/2} \mathbf{M}_k(\mathbf{V}^{\otimes 2})^{-1/2} \mathbf{J}_k^\perp \mathrm{E}_0\left[\text{vec}\left[\sum_{i=1}^n K_{f_1}^{(\ell)}(\tilde{G}_{1k}(d_i^n/\sigma)) \mathbf{U}_i^n \mathbf{U}_i^{n\prime}\right]\right]$$

$$+ \mathcal{J}_k^{(\ell)}(f_1; g_1) \Upsilon_k^{-1}(\mathbf{V}) \mathring{\text{vech}}(\mathbf{v}^{(n)}),$$

$$\mathbf{R}_1^{(n;\ell)} := \frac{1}{2} n^{-1/2} \mathbf{M}_k(\mathbf{V}^{\otimes 2})^{-1/2} \mathbf{J}_k^\perp$$

$$\times \text{vec}\left[\sum_{i=1}^n [K_{f_1}(\tilde{G}_{1k}(d_i^0/\sigma)) - K_{f_1}^{(\ell)}(\tilde{G}_{1k}(d_i^0/\sigma))] \mathbf{U}_i^0 \mathbf{U}_i^{0\prime}\right],$$

$$\mathbf{R}_2^{(n;\ell)} := \frac{1}{2} n^{-1/2} \mathbf{M}_k(\mathbf{V}^{\otimes 2})^{-1/2} \mathbf{J}_k^\perp$$

$$\times \text{vec}\left[\sum_{i=1}^n [K_{f_1}(\tilde{G}_{1k}(d_i^n/\sigma)) - K_{f_1}^{(\ell)}(\tilde{G}_{1k}(d_i^n/\sigma))] \mathbf{U}_i^n \mathbf{U}_i^{n\prime}\right]$$

and

$$\mathbf{R}_3^{(n;\ell)} := (\mathcal{J}_k(f_1, g_1) - \mathcal{J}_k^{(\ell)}(f_1; g_1)) \Upsilon_k^{-1}(\mathbf{V}) \mathring{\text{vech}}(\mathbf{v}^{(n)}).$$

We prove that $\mathbf{C}^{(n)} = o_\mathrm{P}(1)$ (thus completing the proof of Proposition A.1) by establishing that $\mathbf{D}_1^{(n;\ell)}$ and $\mathbf{D}_2^{(n;\ell)}$ are $o_\mathrm{P}(1)$ under $\mathrm{P}_{\boldsymbol{\theta},\sigma^2,\mathbf{V};g_1}^{(n)}$, as $n \to \infty$, for fixed $\ell$ and that $\mathbf{R}_1^{(n;\ell)}$, $\mathbf{R}_2^{(n;\ell)}$ and $\mathbf{R}_3^{(n;\ell)}$ are $o_\mathrm{P}(1)$ under the same sequence of hypotheses, as $\ell \to \infty$, uniformly in $n$. For the sake of convenience, these three results are treated separately (Lemmas A.2, A.3 and A.4). $\square$

LEMMA A.2. *For any fixed $\ell$, $\mathrm{E}_0[\|\mathbf{D}_1^{(n;\ell)}\|^2] = o(1)$ as $n \to \infty$.*

LEMMA A.3. *For any fixed $\ell$, $\mathbf{D}_2^{(n;\ell)} = o(1)$ as $n \to \infty$.*

LEMMA A.4. *As $\ell \to \infty$, uniformly in $n$,*

(i) $\mathbf{R}_1^{(n;\ell)}$ *is* $o_{\mathrm{P}}(1)$ *under* $\mathrm{P}_{\boldsymbol{\theta},\sigma^2,\mathbf{V};g_1}^{(n)}$,

(ii) $\mathbf{R}_2^{(n;\ell)}$ *is* $o_{\mathrm{P}}(1)$ *under* $\mathrm{P}_{\boldsymbol{\theta},\sigma^2,\mathbf{V};g_1}^{(n)}$ *for $n$ sufficiently large,*

(iii) $\mathbf{R}_3^{(n;\ell)}$ *is* $o(1)$.

PROOF OF LEMMA A.2. First, note that

$$\mathbf{D}_1^{(n;\ell)} = \frac{1}{2} n^{-1/2} \mathbf{M}_k (\mathbf{V}^{\otimes 2})^{-1/2} \mathbf{J}_k^\perp \sum_{i=1}^n [\mathbf{T}_i^{(n;\ell)} - \mathrm{E}_0[\mathbf{T}_i^{(n;\ell)}]],$$

where $\mathbf{T}_i^{(n;\ell)} := \mathrm{vec}\,[K_{f_1}^{(\ell)}(\tilde{G}_{1k}(d_i^n/\sigma)) \mathbf{U}_i^n \mathbf{U}_i^{n\prime} - K_{f_1}^{(\ell)}(\tilde{G}_{1k}(d_i^0/\sigma)) \mathbf{U}_i^0 \mathbf{U}_i^{0\prime}]$, $i = 1,\ldots,n$, are i.i.d. Writing $\mathrm{Var}_0$ for variances under $\mathrm{P}_{\boldsymbol{\theta},\sigma^2,\mathbf{V};g_1}^{(n)}$, we have

$$\mathrm{E}_0[\|\mathbf{D}_1^{(n;\ell)}\|^2] \leq C n^{-1} \mathrm{E}_0 \left[ \left\| \sum_{i=1}^n [\mathbf{T}_i^{(n;\ell)} - \mathrm{E}_0[\mathbf{T}_i^{(n;\ell)}]] \right\|^2 \right]$$

$$\leq C n^{-1} \mathrm{tr}\left[ \mathrm{Var}_0 \left[ \sum_{i=1}^n [\mathbf{T}_i^{(n;\ell)} - \mathrm{E}_0[\mathbf{T}_i^{(n;\ell)}]] \right] \right]$$

$$= C \,\mathrm{tr}[\mathrm{Var}_0[\mathbf{T}_1^{(n;\ell)}]] \leq C \mathrm{E}_0[\|\mathbf{T}_1^{(n;\ell)}\|^2],$$

and it only remains to be shown that

$$\begin{aligned}(A.5) \quad \mathrm{E}_0[\|\mathbf{T}_1^{(n;\ell)}\|^2] &= \mathrm{E}_0[\|K_{f_1}^{(\ell)}(\tilde{G}_{1k}(d_1^n/\sigma))\mathrm{vec}\,[\mathbf{U}_1^n \mathbf{U}_1^{n\prime}] \\ &\quad - K_{f_1}^{(\ell)}(\tilde{G}_{1k}(d_1^0/\sigma))\mathrm{vec}\,[\mathbf{U}_1^0 \mathbf{U}_1^{0\prime}]\|^2] = o(1)\end{aligned}$$

as $n \to \infty$. Noting that $\|\mathrm{vec}\,(\mathbf{u}\mathbf{v}')\| = \|\mathbf{u}\|\|\mathbf{v}\|$, we have

$$\|K_{f_1}^{(\ell)}(\tilde{G}_{1k}(d_1^n/\sigma))\mathrm{vec}\,[\mathbf{U}_1^n \mathbf{U}_1^{n\prime}] - K_{f_1}^{(\ell)}(\tilde{G}_{1k}(d_1^0/\sigma))\mathrm{vec}\,[\mathbf{U}_1^0 \mathbf{U}_1^{0\prime}]\|^2$$

$$\leq 2|K_{f_1}^{(\ell)}(\tilde{G}_{1k}(d_1^n/\sigma)) - K_{f_1}^{(\ell)}(\tilde{G}_{1k}(d_1^0/\sigma))|^2 \|\mathrm{vec}\,[\mathbf{U}_1^n \mathbf{U}_1^{n\prime}]\|^2$$

$$+ 2|K_{f_1}^{(\ell)}(\tilde{G}_{1k}(d_1^0/\sigma))|^2 \|\mathrm{vec}\,[\mathbf{U}_1^n \mathbf{U}_1^{n\prime} - \mathbf{U}_1^0 \mathbf{U}_1^{0\prime}]\|^2$$

$$\leq C|K_{f_1}^{(\ell)}(\tilde{G}_{1k}(d_1^n/\sigma)) - K_{f_1}^{(\ell)}(\tilde{G}_{1k}(d_1^0/\sigma))|^2 + C\|\mathbf{U}_1^n - \mathbf{U}_1^0\|^2,$$

for some constant $C$. Lemma A.1(i) and the continuity of $K_{f_1}^{(\ell)} \circ \tilde{G}_{1k}$ together imply that $K_{f_1}^{(\ell)}(\tilde{G}_{1k}(d_1^n/\sigma)) - K_{f_1}^{(\ell)}(\tilde{G}_{1k}(d_1^0/\sigma)) = o_{\mathrm{P}}(1)$, under $\mathrm{P}_{\boldsymbol{\theta},\sigma^2,\mathbf{V};g_1}^{(n)}$, as $n \to \infty$. Since $K_{f_1}^{(\ell)}$ is bounded, this convergence to zero also holds in quadratic mean. Similarly, using Lemma A.1(ii) and the boundedness of $\mathbf{U}_1^0$ and $\mathbf{U}_1^n$, we obtain that $\|\mathbf{U}_1^n - \mathbf{U}_1^0\|$ is $o(1)$ in quadratic mean, as $n \to \infty$, under $\mathrm{P}_{\boldsymbol{\theta},\sigma^2,\mathbf{V};g_1}^{(n)}$. The convergence in (A.5) then follows. □

PROOF OF LEMMA A.3. Letting

$$\mathbf{B}_1^{(n;\ell)} := \frac{1}{2} n^{-1/2} \mathbf{M}_k (\mathbf{V}^{\otimes 2})^{-1/2} \mathbf{J}_k^\perp \mathrm{vec}\left[\sum_{i=1}^n K_{f_1}^{(\ell)}(\tilde{G}_{1k}(d_i^0/\sigma)) \mathbf{U}_i^0 \mathbf{U}_i^{0\prime}\right],$$

one can show that, under $\mathrm{P}_{\boldsymbol{\theta},\sigma^2,\mathbf{V};g_1}^{(n)}$, as $n \to \infty$,

(A.6) $$\mathbf{B}_1^{(n;\ell)} \xrightarrow{\mathcal{L}} \mathcal{N}(\mathbf{0}, \mathrm{E}[(K_{f_1}^{(\ell)}(U))^2] \boldsymbol{\Upsilon}_k^{-1}(\mathbf{V}))$$

[throughout, $U$ stands for a random variable uniformly distributed over $(0,1)$]. Under the sequence of local alternatives $\mathrm{P}_{\boldsymbol{\theta}^n,\sigma^2,\mathbf{V}^n;g_1}^{(n)}$, as $n \to \infty$,

$$\mathbf{B}_1^{(n;\ell)} - \mathcal{J}_k^{(\ell)}(f_1;g_1) \boldsymbol{\Upsilon}_k^{-1}(\mathbf{V}) \overset{\circ}{\mathrm{vech}}(\mathbf{v}^{(n)}) \xrightarrow{\mathcal{L}} \mathcal{N}(\mathbf{0}, \mathrm{E}[(K_{f_1}^{(\ell)}(U))^2] \boldsymbol{\Upsilon}_k^{-1}(\mathbf{V})).$$

Defining $\mathbf{B}_2^{(n;\ell)} := \frac{1}{2} n^{-1/2} \mathbf{M}_k(\mathbf{V}^{\otimes 2})^{-1/2} \mathbf{J}_k^\perp \mathrm{vec}\left[\sum_{i=1}^n K_{f_1}^{(\ell)}(\tilde{G}_{1k}(d_i^n/\sigma)) \mathbf{U}_i^n \mathbf{U}_i^{n\prime}\right]$, it follows from ULAN that, under $\mathrm{P}_{\boldsymbol{\theta},\sigma^2,\mathbf{V};g_1}^{(n)}$, as $n \to \infty$,

(A.7)
$$\mathbf{B}_2^{(n;\ell)} + \mathcal{J}_k^{(\ell)}(f_1;g_1) \boldsymbol{\Upsilon}_k^{-1}(\mathbf{V}) \overset{\circ}{\mathrm{vech}}(\mathbf{v}^{(n)})$$
$$\xrightarrow{\mathcal{L}} \mathcal{N}(\mathbf{0}, \mathrm{E}[(K_{f_1}^{(\ell)}(U))^2] \boldsymbol{\Upsilon}_k^{-1}(\mathbf{V})).$$

Now, from (A.6) and the fact that, under $\mathrm{P}_{\boldsymbol{\theta},\sigma^2,\mathbf{V};g_1}^{(n)}$, $\mathbf{D}_1^{(n;\ell)} = \mathbf{B}_2^{(n;\ell)} - \mathbf{B}_1^{(n;\ell)} - \mathrm{E}_0[\mathbf{B}_2^{(n;\ell)}] = o_\mathrm{P}(1)$ as $n \to \infty$ (Lemma A.2), we obtain that

(A.8) $$\mathbf{B}_2^{(n;\ell)} - \mathrm{E}_0[\mathbf{B}_2^{(n;\ell)}] \xrightarrow{\mathcal{L}} \mathcal{N}(\mathbf{0}, \mathrm{E}[(K_{f_1}^{(\ell)}(U))^2] \boldsymbol{\Upsilon}_k^{-1}(\mathbf{V}))$$

as $n \to \infty$ under $\mathrm{P}_{\boldsymbol{\theta},\sigma^2,\mathbf{V};g_1}^{(n)}$. Comparing (A.7) and (A.8), it follows that $\mathbf{D}_2^{(n;\ell)} = \mathrm{E}_0[\mathbf{B}_2^{(n;\ell)}] + \mathcal{J}_k^{(\ell)}(f_1;g_1) \boldsymbol{\Upsilon}_k^{-1}(\mathbf{V}) \overset{\circ}{\mathrm{vech}}(\mathbf{v}^{(n)})$ is $o(1)$ as $n \to \infty$. □

We now complete the proof of Proposition A.1 by proving Lemma A.4.

PROOF OF LEMMA A.4. (i) In view of the independence, under $\mathrm{P}_{\boldsymbol{\theta},\sigma^2,\mathbf{V};g_1}^{(n)}$, between the $d_i^0$'s and the $\mathbf{U}_i^0$'s, we obtain, for all $n$,

$$\mathrm{E}_0[\|\mathbf{R}_1^{(n;\ell)}\|^2] \leq \frac{C}{n} \sum_{i=1}^n \mathrm{E}_0[[K_{f_1}(\tilde{G}_{1k}(d_i^0/\sigma)) - K_{f_1}^{(\ell)}(\tilde{G}_{1k}(d_i^0/\sigma))]^2]$$

$$\times \mathrm{E}_0[[\mathrm{vec}\,\mathbf{U}_i^0 \mathbf{U}_i^{0\prime}]' \mathbf{J}_k^\perp [\mathrm{vec}\,\mathbf{U}_i^0 \mathbf{U}_i^{0\prime}]]$$

(A.9)
$$= \frac{C(k-1)}{kn} \sum_{i=1}^n \mathrm{E}_0[[K_{f_1}(\tilde{G}_{1k}(d_i^0/\sigma)) - K_{f_1}^{(\ell)}(\tilde{G}_{1k}(d_i^0/\sigma))]^2]$$

$$= \frac{C(k-1)}{k} \int_0^1 [K_{f_1}(u) - K_{f_1}^{(\ell)}(u)]^2 \, du.$$

Now, $K_{f_1}^{(\ell)}(u)$ converges to $K_{f_1}(u)$ for all $u \in (0,1)$. Also, since $|K_{f_1}^{(\ell)}(u)|$ is bounded by $|K_{f_1}(u)|$, for all $\ell \geq L$, the integrand in (A.9) is bounded (uniformly in $\ell$) by $4|K_{f_1}(u)|^2$ which is integrable on $(0,1)$. The Lebesgue dominated convergence theorem thus yields that $\mathrm{E}_0[\|\mathbf{R}_1^{(n;\ell)}\|^2] = o(1)$ as $\ell \to \infty$. This convergence is uniform in $n$, since the constant $C$ in (A.9) does not depend on $n$.

(ii) The claim in part (ii) of the lemma is the same as in part (i), except that $d_i^n$ and $\mathbf{U}_i^n$ replace $d_i^0$ and $\mathbf{U}_i^0$, respectively. Accordingly, part (ii) holds under $\mathrm{P}_{\boldsymbol{\theta}^n,\sigma^2,\mathbf{V}^n;g_1}^{(n)}$. That it also holds under $\mathrm{P}_{\boldsymbol{\theta},\sigma^2,\mathbf{V};g_1}^{(n)}$ follows from Lemma 3.5 of [23].

(iii) Note that $|\mathcal{J}_k(f_1,g_1) - \mathcal{J}_k^{(\ell)}(f_1;g_1)|^2 = |\int_0^1 (K_{f_1}(u) - K_{f_1}^{(\ell)}(u))K_{g_1}(u)\,du|^2$
$\leq \mathcal{J}_k(g_1)\int_0^1 |K_{f_1}(u) - K_{f_1}^{(\ell)}(u)|^2\,du$. Again, $|K_{f_1}^{(\ell)}(u) - K_{f_1}(u)|^2 \leq 4|K_{f_1}(u)|^2$ with $\int_0^1 |K_{f_1}(u)|^2\,du < \infty$. Pointwise convergence of $(K_{f_1}^{(\ell)})$ to $K$ implies that $\mathcal{J}_k(f_1,g_1) - \mathcal{J}_k^{(\ell)}(f_1;g_1) = o(1)$ as $\ell \to \infty$. The result then follows from the boundedness of $(\mathbf{v}^{(n)})$. □

### A.2. Properties of $R$-estimators.

PROOF OF PROPOSITION 3.1. (i) The asymptotic representations (3.5) and (3.6) are just restatements of (3.2) and (3.3), to which we refer for the proof. The convergence in (3.7) then readily results from part (iii) of Proposition 2.1. As for (3.8), it follows directly from the fact that $\mathrm{vec}(\underset{\sim}{\mathbf{V}}_{f_1\#}^{(n)} - \mathbf{V}) = \mathbf{M}_k' \,\mathrm{vech}(\underset{\sim}{\mathbf{V}}_{f_1\#}^{(n)} - \mathbf{V})$ and the definition of $\mathbf{Q}_k(\mathbf{V})$.

(ii) Semiparametric efficiency follows from the fact that $\mathcal{J}_k(f_1,f_1) = \mathcal{J}_k(f_1)$ so that under $\mathrm{P}_{\sigma^2,\mathbf{V};f_1}^{(n)}$, the asymptotic variance in (3.7) reduces to $\mathcal{J}_k(f_1)^{-1}\boldsymbol{\Upsilon}_k(\mathbf{V})$, the inverse of the efficient information matrix $\boldsymbol{\Gamma}_{f_1}^*(\mathbf{V})$.

(iii) From (3.4) and (3.1), [with $R_i = R_i^{(n)}(\mathbf{V}_\#^{(n)})$ and $\mathbf{U}_i = \mathbf{U}_i^{(n)}(\mathbf{V}_\#^{(n)})$],

$$\mathrm{vech}(\underset{\sim}{\mathbf{V}}_{f_1\#}^{(n)}) = \mathrm{vech}(\mathbf{V}_\#^{(n)}) + \frac{k(k+2)}{n^{1/2}\mathcal{J}_k(f_1,g_1)}\mathbf{N}_k\mathbf{Q}_k(\mathbf{V}_\#^{(n)})\mathbf{N}_k'\underset{\sim}{\boldsymbol{\Delta}}_{f_1}^{(n)}(\mathbf{V}_\#^{(n)})$$

$$= \mathrm{vech}(\mathbf{V}_\#^{(n)}) + \frac{k(k+2)}{2n\mathcal{J}_k(f_1,g_1)}\mathbf{N}_k\mathbf{Q}_k(\mathbf{V}_\#^{(n)})((\mathbf{V}_\#^{(n)})^{\otimes 2})^{-1/2}$$

$$\times \sum_{i=1}^n \left[K_{f_1}\left(\frac{R_i}{n+1}\right)\mathrm{vec}(\mathbf{U}_i\mathbf{U}_i') - \frac{m_{f_1}^{(n)}}{k}\mathrm{vec}(\mathbf{I}_k)\right],$$

where we used the fact that (see Section 4.2 for the definition of $\mathbf{e}_{k^2,1}$)

$$\mathbf{Q}_k(\mathbf{V})\mathbf{N}_k'\mathbf{M}_k = \mathbf{Q}_k(\mathbf{V})$$

$$= [\mathbf{I}_{k^2} - (\operatorname{vec}\mathbf{V})\mathbf{e}'_{k^2,1}][\mathbf{I}_{k^2} + \mathbf{K}_k](\mathbf{V}^{\otimes 2})[\mathbf{I}_{k^2} - (\operatorname{vec}\mathbf{V})\mathbf{e}'_{k^2,1}]'$$

$$= [\mathbf{I}_{k^2} + \mathbf{K}_k](\mathbf{V}^{\otimes 2}) - 2(\mathbf{V}^{\otimes 2})\mathbf{e}_{k^2,1}(\operatorname{vec}\mathbf{V})'$$

$$\quad - 2(\operatorname{vec}\mathbf{V})\mathbf{e}'_{k^2,1}(\mathbf{V}^{\otimes 2}) + 2(\operatorname{vec}\mathbf{V})(\operatorname{vec}\mathbf{V})';$$

see the proof of Lemma HP3.1. Routine algebra yields

$$\operatorname{vech}(\widehat{\overset{\circ}{\underset{\sim}{\mathbf{V}}}}{}^{(n)}_{f_1}) = \operatorname{vech}(\mathbf{V}^{(n)}_{\#}) + \frac{k(k+2)}{\alpha^*}\mathbf{N}_k[\mathbf{I}_{k^2} - (\operatorname{vec}\mathbf{V}^{(n)}_{\#})\mathbf{e}'_{k^2,1}]((\mathbf{V}^{(n)}_{\#})^{\otimes 2})^{1/2}$$

$$\times \left(\frac{1}{n}\sum_{i=1}^n K_{f_1}\left(\frac{R_i}{n+1}\right)\operatorname{vec}(\mathbf{U}_i\mathbf{U}'_i)\right)$$

(A.10)

$$= \operatorname{vech}(\mathbf{V}^{(n)}_{\#}) + \frac{k(k+2)}{\alpha^*}\mathbf{N}_k[\mathbf{I}_{k^2} - (\operatorname{vec}\mathbf{V}^{(n)}_{\#})\mathbf{e}'_{k^2,1}]\operatorname{vec}(\underset{\sim}{\mathbf{W}}^{(n)}_{f_1\#})$$

$$= \operatorname{vech}(\mathbf{V}^{(n)}_{\#}) + \frac{k(k+2)}{\alpha^*}\mathbf{N}_k\operatorname{vec}(\underset{\sim}{\mathbf{W}}^{(n)}_{f_1\#} - (\underset{\sim}{\mathbf{W}}^{(n)}_{f_1\#})_{11}\mathbf{V}^{(n)}_{\#}),$$

which establishes the result since $\operatorname{vech}\mathbf{v} = \operatorname{vech}\mathbf{w}$ if and only if $\mathbf{v} = \mathbf{w}$ for all $k \times k$ symmetric matrices $\mathbf{v} = (v_{ij})$, $\mathbf{w} = (w_{ij})$ such that $v_{11} = w_{11}$.

(iv) Due to the identification constraints, the population covariance matrix under $P^{(n)}_{\sigma^2,\mathbf{V};g_1}$ with finite second-order moments is not $\boldsymbol{\Sigma} := \sigma^2\mathbf{V}$, but $\eta\boldsymbol{\Sigma} := k^{-1}\sigma^2 D_k(g_1)\mathbf{V}$. Provided that $\kappa_k(g_1) < \infty$, the multivariate central limit theorem yields $n^{1/2}\operatorname{vec}(\boldsymbol{\Sigma}^{(n)} - \eta\boldsymbol{\Sigma}) \overset{\mathcal{L}}{\to} \mathcal{N}(\mathbf{0}, \mathbf{A})$, where

$$\mathbf{A} := \frac{\sigma^4 E_k(g_1)}{k(k+2)}[\mathbf{I}_{k^2} + \mathbf{K}_k](\mathbf{V}^{\otimes 2}) + \frac{\sigma^4 \kappa_k(g_1) D_k^2(g_1)}{k^2}(\operatorname{vec}\mathbf{V})(\operatorname{vec}\mathbf{V})'.$$

Now, applying Slutsky's lemma, we obtain, as $n \to \infty$, under $P^{(n)}_{\sigma^2,\mathbf{V};g_1}$,

$$n^{1/2}\operatorname{vec}(\mathbf{V}^{(n)}_{\mathcal{G}} - \mathbf{V}) = \frac{1}{\eta\Sigma_{11}}[\mathbf{I}_{k^2} - (\operatorname{vec}\mathbf{V})\mathbf{e}'_{k^2,1}][n^{1/2}\operatorname{vec}(\boldsymbol{\Sigma}^{(n)} - \eta\boldsymbol{\Sigma})] + o_P(1)$$

$$\overset{\mathcal{L}}{\to} \mathcal{N}\left(\mathbf{0}, \frac{1}{\eta^2\sigma^4}[\mathbf{I}_{k^2} - (\operatorname{vec}\mathbf{V})\mathbf{e}'_{k^2,1}]\mathbf{A}[\mathbf{I}_{k^2} - (\operatorname{vec}\mathbf{V})\mathbf{e}'_{k^2,1}]'\right),$$

where the covariance matrix, after lengthy but standard algebra, reduces to $(1 + \kappa_k(g_1))\mathbf{Q}_k(\mathbf{V})$, yielding the desired result; see also [36].

(v) The asymptotic covariance matrices of $\operatorname{vec}(\underset{\sim}{\mathbf{V}}^{(n)}_{f_1\#})$ in (3.8) and $\operatorname{vec}(\mathbf{V}^{(n)}_{\mathcal{G}})$ in (iv) are proportional; ARE's with respect to $\mathbf{V}^{(n)}_{\mathcal{G}}$ in (v) follow directly as ratios of the corresponding proportionality factors. ARE's with respect to $\mathbf{V}^{(n)}_T$ follow from the fact that in the normalization adopted [i.e., $(\mathbf{V}^{(n)}_T)_{11} = 1$], $n^{1/2}\operatorname{vec}(\mathbf{V}^{(n)}_T - \mathbf{V})$ is asymptotically normal with mean zero and covariance matrix $((k+2)/k)\mathbf{Q}_k(\mathbf{V})$. □

### A.3. Asymptotic equivariance.

PROOF OF PROPOSITION 5.1.
(i) We first prove that

(A.11) $$\underset{\sim}{\mathbf{W}}_{f_1\#}^{(n)} - \mathbf{V}_{\#}^{(n)} = O_{\mathrm{P}}(n^{-1/2}),$$

under $\mathcal{P}^{(n)}$, as $n \to \infty$ [recall that $\underset{\sim}{\mathbf{W}}_{f_1\#}^{(n)} := \underset{\sim}{\mathbf{W}}_{f_1}^{(n)}(\mathbf{V}_{\#}^{(n)})$]. To this end, define

$$\underset{\sim}{\mathbf{T}}_{f_1}^{(n)}(\mathbf{V}) := n^{-1/2}(\mathbf{V}^{\otimes 2})^{1/2} \sum_{i=1}^{n} \left[ K_{f_1}\left(\frac{R_i}{n+1}\right) \mathrm{vec}(\mathbf{U}_i\mathbf{U}_i') - \frac{m_{f_1}^{(n)}}{k} \mathrm{vec}(\mathbf{I}_k) \right]$$

[with $R_i = R_i^{(n)}(\mathbf{V})$ and $\mathbf{U}_i = \mathbf{U}_i^{(n)}(\mathbf{V})$], which is asymptotically normal with mean zero and covariance matrix $\mathcal{J}_k(f_1)\mathbf{H}_k(\mathbf{V})$, where

$$\mathbf{H}_k(\mathbf{V}) := \frac{1}{k(k+2)}(\mathbf{V}^{\otimes 2})^{1/2}\left[\mathbf{I}_{k^2} + \mathbf{K}_k - \frac{2}{k}\mathbf{J}_k\right](\mathbf{V}^{\otimes 2})^{1/2}.$$

Proceeding exactly as in the proof of Proposition A.1, we obtain that for any bounded sequence $\mathbf{v}^{(n)}$ of symmetric matrices such that $v_{11}^{(n)} = 0$, the difference

$$\underset{\sim}{\mathbf{T}}_{f_1}^{(n)}(\mathbf{V} + n^{-1/2}\mathbf{v}^{(n)}) - \underset{\sim}{\mathbf{T}}_{f_1}^{(n)}(\mathbf{V}) + \frac{1}{2}\mathcal{J}_k(f_1,g_1)\mathbf{H}_k(\mathbf{V})(\mathbf{V}^{\otimes 2})^{-1}\mathrm{vec}(\mathbf{v}^{(n)})$$

is $o_{\mathrm{P}}(1)$, under $\mathrm{P}_{\sigma^2,\mathbf{V};g_1}^{(n)}$, as $n \to \infty$. The local discreteness of $\mathbf{V}_{\#}^{(n)}$ allows us to replace the nonrandom quantity $\mathbf{V}^{(n)} = \mathbf{V} + n^{-1/2}\mathbf{v}^{(n)}$ with the random one $\mathbf{V}_{\#}^{(n)}$ (see, e.g., [30], Lemma 4.4), yielding

$$\underset{\sim}{\mathbf{T}}_{f_1}^{(n)}(\mathbf{V}_{\#}^{(n)}) - \underset{\sim}{\mathbf{T}}_{f_1}^{(n)}(\mathbf{V}) + \frac{1}{2}\mathcal{J}_k(f_1,g_1)\mathbf{H}_k(\mathbf{V})(\mathbf{V}^{\otimes 2})^{-1}n^{1/2}\mathrm{vec}(\mathbf{V}_{\#}^{(n)} - \mathbf{V}) = o_{\mathrm{P}}(1),$$

under $\mathrm{P}_{\sigma^2,\mathbf{V};g_1}^{(n)}$, as $n \to \infty$. This establishes (A.11) since

$$n^{1/2}\mathrm{vec}(\underset{\sim}{\mathbf{W}}_{f_1\#}^{(n)} - \mathbf{V}_{\#}^{(n)}) = \underset{\sim}{\mathbf{T}}_{f_1}^{(n)}(\mathbf{V}_{\#}^{(n)}) + n^{1/2}k^{-1}(m_{f_1}^{(n)} - k)\mathrm{vec}(\mathbf{V}_{\#}^{(n)})$$

$$= n^{1/2}k^{-1}(m_{f_1}^{(n)} - k)\mathrm{vec}(\mathbf{V}_{\#}^{(n)}) + \underset{\sim}{\mathbf{T}}_{f_1}^{(n)}(\mathbf{V})$$

(A.12)
$$- \frac{1}{2}\mathcal{J}_k(f_1,g_1)\mathbf{H}_k(\mathbf{V})(\mathbf{V}^{\otimes 2})^{-1}n^{1/2}\mathrm{vec}(\mathbf{V}_{\#}^{(n)} - \mathbf{V})$$

$$+ o_{\mathrm{P}}(1)$$

(still under $\mathrm{P}_{\sigma^2,\mathbf{V};g_1}^{(n)}$, as $n \to \infty$) and since the square-integrability of $K_{f_1}$ over $(0,1)$ implies that $m_{f_1}^{(n)} - k = m_{f_1}^{(n)} - \int_0^1 K_{f_1}(u)\,du = o(n^{-1/2})$ (see the proof of Proposition 3.2(i) in [16]).

Now, denoting by $\underset{\sim}{\mathbf{V}}_{f_1\#}^{(n)} := \underset{\sim}{\mathbf{V}}_{f_1\#}^{(n)}(c_0)$ the pseudo-estimator defined in (3.1), it follows from (A.11) that [letting $b := k(k+2)\mathcal{J}_k^{-1}(f_1, g_1)$]

$$\text{vec}(\underset{\approx}{\mathbf{V}}_{f_1\#}^{(n)} - \underset{\sim}{\mathbf{V}}_{f_1\#}^{(n)}) = (-b^2(\underset{\sim}{\mathbf{W}}_{f_1\#}^{(n)} - \mathbf{V}_{\#}^{(n)})_{11})(1 + b(\underset{\sim}{\mathbf{W}}_{f_1\#}^{(n)} - \mathbf{V}_{\#}^{(n)})_{11})^{-1}$$

$$\times [\mathbf{I}_{k^2} - (\text{vec}\,\mathbf{V}_{\#}^{(n)})\mathbf{e}'_{k^2,1}]\text{vec}(\underset{\sim}{\mathbf{W}}_{f_1\#}^{(n)} - \mathbf{V}_{\#}^{(n)})$$

is $o_P(n^{-1/2})$, under $\mathcal{P}^{(n)}$, as $n \to \infty$. This yields the result since, in Section 4.2, we proved that $\underset{\sim}{\mathbf{V}}_{f_1\#}^{(n)} - \underset{\sim}{\widehat{\mathbf{V}}}_{f_1\#}^{(n)} = o_P(n^{-1/2})$, under $\mathcal{P}^{(n)}$, as $n \to \infty$.

(ii) If $\mathbf{V}^{(n)}$ is strictly affine-equivariant [in the sense of (5.1)], then using the same notation as in Section 5, $(\mathbf{V}^{(n)}(\mathbf{M}, \mathbf{a}))^{1/2} = d\mathbf{M}(\mathbf{V}^{(n)})^{1/2}\mathbf{O}$ for some $d > 0$ and some $k \times k$ orthogonal matrix $\mathbf{O}$ (see, e.g., [40]). The strict affine-equivariance of the practical implementation $\underset{\approx}{\mathbf{V}}_{f_1}^{(n)} = \lim_{m\to\infty} \underset{\approx}{\mathbf{V}}_{f_1\#}^{(n)}(c_{0,m})$ [which is based on $\mathbf{V}_T^{(n)}$ and $\underset{\sim}{\mathbf{W}}_{f_1}^{(n)}(\mathbf{V}_T^{(n)})$ instead of $\mathbf{V}_{\#}^{(n)}$ and $\underset{\sim}{\mathbf{W}}_{f_1\#}^{(n)}$] follows. $\square$

Département de Mathématique
Institute for Research in Statistics and E.C.A.R.E.S.
Université Libre de Bruxelles
Campus de la Plaine CP 210
B-1050 Bruxelles
Belgium
E-mail: mhallin@ulb.ac.be

Tampere School of Public Health
University of Tampere
FIN-33014 Tampere
Finland
E-mail: Hannu.Oja@uta.fi

Département de Mathématique
Institute for Research in Statistics and E.C.A.R.E.S.
Université Libre de Bruxelles
Campus de la Plaine CP 210
B-1050 Bruxelles
Belgium
E-mail: dpaindav@ulb.ac.be
URL: http://homepages.ulb.ac.be/~dpaindav